\documentclass{amsart}

\usepackage{microtype}
\usepackage{amsmath}
\usepackage{amssymb}
\usepackage{mathtools}
\usepackage{bm}
\usepackage{enumerate}
\usepackage{eqlist}
\usepackage{enumitem}
\usepackage{soul}

\usepackage[dvipsnames]{xcolor}
\usepackage[most]{tcolorbox}
\usepackage{./paragraphs}
\usepackage[colorlinks,unicode,psdextra]{hyperref}
\hypersetup{citecolor = {Green}}
\hypersetup{linkcolor = {Red}}
\usepackage[margin=2\parindent]{caption}

\usepackage{aliascnt}

\usepackage{catoptions}
\makeatletter

\def\Autoref#1{%
  \begingroup
  \edef\reserved@a{\cpttrimspaces{#1}}%
  \ifcsndefTF{r@#1}{%
    \xaftercsname{\expandafter\testreftype\@fourthoffive}
      {r@\reserved@a}.\\{#1}%
  }{%
    \ref{#1}%
  }%
  \endgroup
}
\def\testreftype#1.#2\\#3{%
  \ifcsndefTF{#1autorefname}{%
    \def\reserved@a##1##2\@nil{%
      \uppercase{\def\ref@name{##1}}%
      \csn@edef{#1autorefname}{\ref@name##2}%
      \autoref{#3}%
    }%
    \reserved@a#1\@nil
  }{%
    \autoref{#3}%
  }%
}
\makeatother
\newcommand*{\Aref}[1]{\Autoref{#1}}

\DeclareCaptionLabelFormat{nospace}{#1#2}
\let\thefig\thefigure
\renewcommand{\thefigure}{Figure \textup{\text{\thefig}}}

\numberwithin{equation}{section}
\makeatletter
\newcommand{\leqnos}{\tagsleft@true\let\veqno\@@leqno}
\newcommand{\reqnos}{\tagsleft@false\let\veqno\@@eqno}
\leqnos
\makeatother

\renewcommand{\textsc}[1]{{\fontseries{m}\fontshape{sc}\selectfont #1}}
\renewcommand{\textbfsc}[1]{{\fontseries{m}\fontshape{sc}\selectfont #1}}

\usepackage{amsrefs}
\newcommand{\biblio}{%
\bibliography{./sources}
}

\allowdisplaybreaks
\newcommand{\ubnote}[3]{%
    \underbrace{#1}_{\clap{\parbox{#2}{\vspace{1mm}\footnotesize\centering #3}}}%
}
\makeatletter
\long\def\@firstofthree#1#2#3{#1}
\newcommand{\smashub}[2]{%
  \smash{#1}%
  \vphantom{\let\ubnote\@firstofthree#1}%
  \gdef\ubfixvsp{\vphantom{#1}}%
}
\makeatother
\newcommand{\subnote}[3]{\smashub{\ubnote{#1}{#2}{#3}}{\subnotefixvsp}}

\newcommand{\mnote}[3][t]{%
    \parbox[#1]{#2}{\raggedright\footnotesize #3}
}
\newcommand*{\R}{\mathbf{R}}
\renewcommand*{\C}{\mathbf{C}}

\newcommand*{\Z}{\mathbf{Z}}
\newcommand*{\Q}{\mathbf{Q}}
\newcommand*{\e}{\mathrm{e}}
\newcommand*{\im}{\mathfrak{i}}
\renewcommand*{\epsilon}{\varepsilon}
\renewcommand*{\phi}{\varphi}
\newcommand*{\di}{\mathrm{d}}
\newcommand*{\f}{\frac}
\newcommand*{\defeq}{\mathrel{\vcenter{\baselineskip0.55ex \lineskiplimit0pt\hbox{\scriptsize.}\hbox{\scriptsize.}}{=}}}
\newcommand*{\eqdef}{\mathrel{{=}\vcenter{\baselineskip0.55ex \lineskiplimit0pt\hbox{\scriptsize.}\hbox{\scriptsize.}}}}
\newcommand*{\fdtr}{\longrightarrow}
\newcommand*{\Mapsto}{\longmapsto}
\DeclareMathOperator{\sgn}{sgn}
\newcommand{\abs}[1]{\left| #1 \right|}
\newcommand{\p}[1]{\left( #1 \right)}
\newcommand{\cur}[1]{\left\{ #1 \right\}}
\newcommand{\s}[1]{\left[ #1 \right]}
\newcommand{\ioo}[1]{\left( #1 \right)}

\newcommand{\ol}[1]{\overline{#1}}
\newcommand*{\msp}{\enskip}
\DeclareMathOperator{\sech}{sech}
\DeclareMathOperator{\arccosh}{arccosh}

\DeclareMathOperator{\arcsec}{arcsec}

\newcommand{\ff}[1]{\mathop{\bm{f}_{#1}}}

\newcommand{\ee}[1]{\mathop{\bm{e}_{#1}}}
\newcommand{\eeh}[1]{\mathop{\hat{\bm{e}}_{#1}}}
\newcommand{\eec}[1]{\mathop{\check{\bm{e}}_{#1}}}
\newcommand{\ees}[1]{\mathop{\bm{e}_{#1}^{*}}}

\newcommand{\Sph}{\mathbf{S}}

\newcommand{\mbinom}[2]{\scalebox{0.8}{$\dbinom{#1}{#2}$}}
\newcommand*{\mpn}[1]{\mathrlap{\quad#1}}

\begin{document}
\setcounter{section}{-1}

\author[M.R.~Jimenez]{Michael Robert Jimenez}
\address{Technische Universit\"{a}t Wien}
\email{michael.jimenez@tuwien.ac.at}

\title[Surfaces with a Relation between their Curvature Radii]{Note on Surfaces of Revolution\\ with an Affine-Linear Relation\\ between their Curvature Radii}

\date{}

\begin{abstract}
  This note derives parametrizations for surfaces of revolution that satisfy an affine-linear relation between their respective curvature radii.  Alongside, parametrizations for the uniform normal offsets of those surfaces are obtained.  Those parametrizations are found explicitly for a countably-infinite many of them, and of those, it is shown which are algebraic.  Lastly, for those surfaces which have a constant ratio of principal curvatures, parametrizations with a constant angle between the parameter curves are found.
\end{abstract}

\maketitle

\section{Introduction}

In one of the foundational papers \cite{Hopf} on the now-named Weingarten surfaces, H.~Hopf found an integral formula for surfaces of revolution which have a constant ratio of principal curvatures, cf.\ equation (29) ibidem, and discussion in English in example 3.27 of \cite{KuehnelDiffGeom}.  For the particular case when the ratio of the principal curvatures is a positive constant, these surfaces are called ``Mylar balloons'', which are discussed in \cite{Hadzhilazova-Mladenov,Mladenov-Opera_balloon,Mladenov-Opera_deform,Pulov-Hadzhilazova-Mladenov}.  And, when that ratio is a negative constant, generalized Chebyshev nets are found in \cite{Staeckel,Riveros-CorroI,Riveros-CorroII}, dropping the surface-of-revolution qualification.  There is also a relation between surfaces that have a constant ratio of principal curvatures with those satisfying equality between the curvature of the second fundamental form and the mean curvature, $K_{\mathrm{II}} = H$, as discussed in \cite{Kuehnel,Kuehnel-Steller} for surfaces of revolution, and in \cite{Baikoussis-Koufogiorgos} for helical surfaces.  Interestingly, surfaces with a constant ratio of principal curvatures appear in solutions to the non-linear Schr\"{o}dinger equation and modified Korteweg-de-Vries equation \cite{Ceyhan-Fokas-Guerses}.

The purpose of this note is to obtain another integral formula for those surfaces of revolution Hopf considered. The objective being therewith, that those surfaces are easier to explicitly compute, with the added benefit of being able to directly obtain their uniform normal offsets.  In obtaining the explicit parametrizations, it will be possible to prove that some of these surfaces are algebraic, while some others are transcendental.  Moreover, the parametrizations with a constant angle between the parameter curves will be found, which for the negatively-curved surfaces will be their asymptotic parametrization.

In particular, this note will focus on surfaces of revolution with a constant ratio of principal curvatures, as part of a larger class of Weingarten surfaces: namely, surfaces that have an affine-linear relationship between their principal curvature radii $\rho_{1}$ and $\rho_{2}$,
\[\rho_{1} + m \rho_{2} = c ,\]
for $m , c \in \R$.  When $c = 0$, this relationship reduces to a constant ratio of principal curvature radii, or thusly, of principal curvatures, while for nearly all $m$, $c \ne 0$ corresponds to a uniform normal offset of the $c = 0$ surface for a given $m$. Surfaces satisfying this relation, will be called of \emph{$(m , c)$-type}.  Other papers that consider classes of Weingarten surfaces containing surfaces that have a constant ratio of principal curvatures are, for example, \cite{Lopez_Htype,Lopez,Lopez-Pampano,Lopez_EandH,Boyacioglu-Lopez,Boyacioglu-Lopez_spacelike}.

For such surfaces when $c = 0$ that are negatively-curved ($m > 0$), their relevance comes from a previous joint paper \cite{Jimenez-Mueller-Pottmann} of the author of this note, wherein the application of these surfaces to architecture is discussed, as fa\c{c}ades with support structures coming in a geometrically-natural way from its asymptotic curves, \`{a} la \cite{Tang-Killian-Bo-Wallner-Pottmann}.  Furthermore, the ``beams'' of those support structures are perpendicular to their fa\c{c}ade surface, which, when they are of uniform height, makes them perpendicular to a uniform normal offset of the fa\c{c}ade surface.

In \Aref{sec:support}, the support parametrization of the profile curve will be constructed from section 6.8 of \cite{Gray}, from which a differential equation in $\rho_{2}$ for these surfaces will be determined.  This differential equation will be solved in \Aref{sec:m-1} for $m = - 1$, as well as include a discussion about that particular subfamily of solution surfaces.  The remaining cases, when $m \ne 0 , - 1$, will be obtained with an integral formula derived in \Aref{sec:mne-1}.  Those integrals will be explicitly computed for $m \in \Z$ in \Aref{sec:alg}, with remarks about which are algebraic, proving a conjecture from \cite{Pottmann}.  Following that, \Aref{sec:param} will find for those surfaces, which have a constant ratio of principal curvatures, parametrizations with a constant angle between the parameter curves, using the manuscript \cite{Pottmann} as a foundation.  Lastly, \Aref{sec:examples} will show some examples.

\section{Support Parametrization}
\label{sec:support}

In this note, a surface of \emph{$(m , c)$-type}, for some constants $m , c \in \R$, will be a surface of revolution $\Sigma$ swept out by the profile curve $\p{r \p{\theta}, h \p{\theta}}$, parametrized as
\[\p{r \p{\theta} \cos \p{t} ,\quad r \p{\theta} \sin \p{t} ,\quad h \p{\theta}} \mpn{,}\]
that satisfies an affine-linear relation between its principal curvature radii $\rho_{1}$, $\rho_{2}$,
\[\rho_{1} + m \rho_{2} = c \mpn{.}\label{eq:aff-lin}\tag{$*$}\]
Without loss of generality, $\Sigma$ will be so oriented that $\rho_{1}$ corresponds to the curvature radius of the profile curve, and $\rho_{2}$ is the curvature of the parallel circles.  The case when $m = 0$ is excluded, as it will be considered singular: $\Sigma$ is swept out by a profile curve with constant curvature as a plane curve, which is a circle.  Also, when $m \ne 0$ and either $\rho_{1}$ or $\rho_{2}$ is constant, then they both must be constant, and as a consequence of corresponding to a surface of revolution, they must also be equal, forcing the resulting surface to be a sphere, with $m = -1$ and $c = 0$.

\begin{figure}[h!]
  \centering{%
    \begin{minipage}[m]{0.9\textwidth}
      \includegraphics[width=\textwidth]{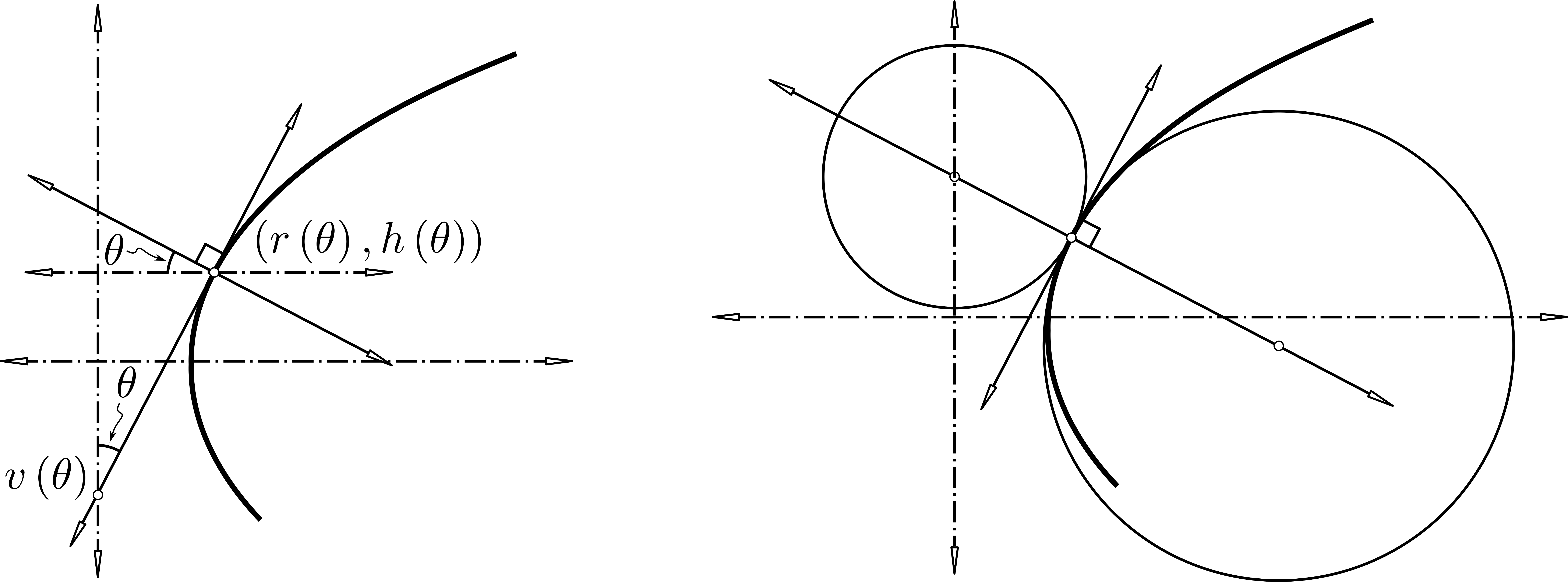}
    \end{minipage}
  }
  \caption{%
    Depictions of a segment of a profile curve $\p{r \p{\theta}, h \p{\theta}}$ of a surface of revolution $\Sigma$, whose axis of revolution is the depicted vertical axis.
    \emph{Left}: A member of the one-parameter family of lines, that are enveloped by that profile curve.
    \emph{Right}: A cross section of the principal curvature spheres of $\Sigma$, with respect to that profile curve.
  }
   \label{fig:radii}
\end{figure}

\begin{figure}[h!]
  \centering{%
    \begin{minipage}[m]{0.3\textwidth}
      \includegraphics[width=\textwidth]{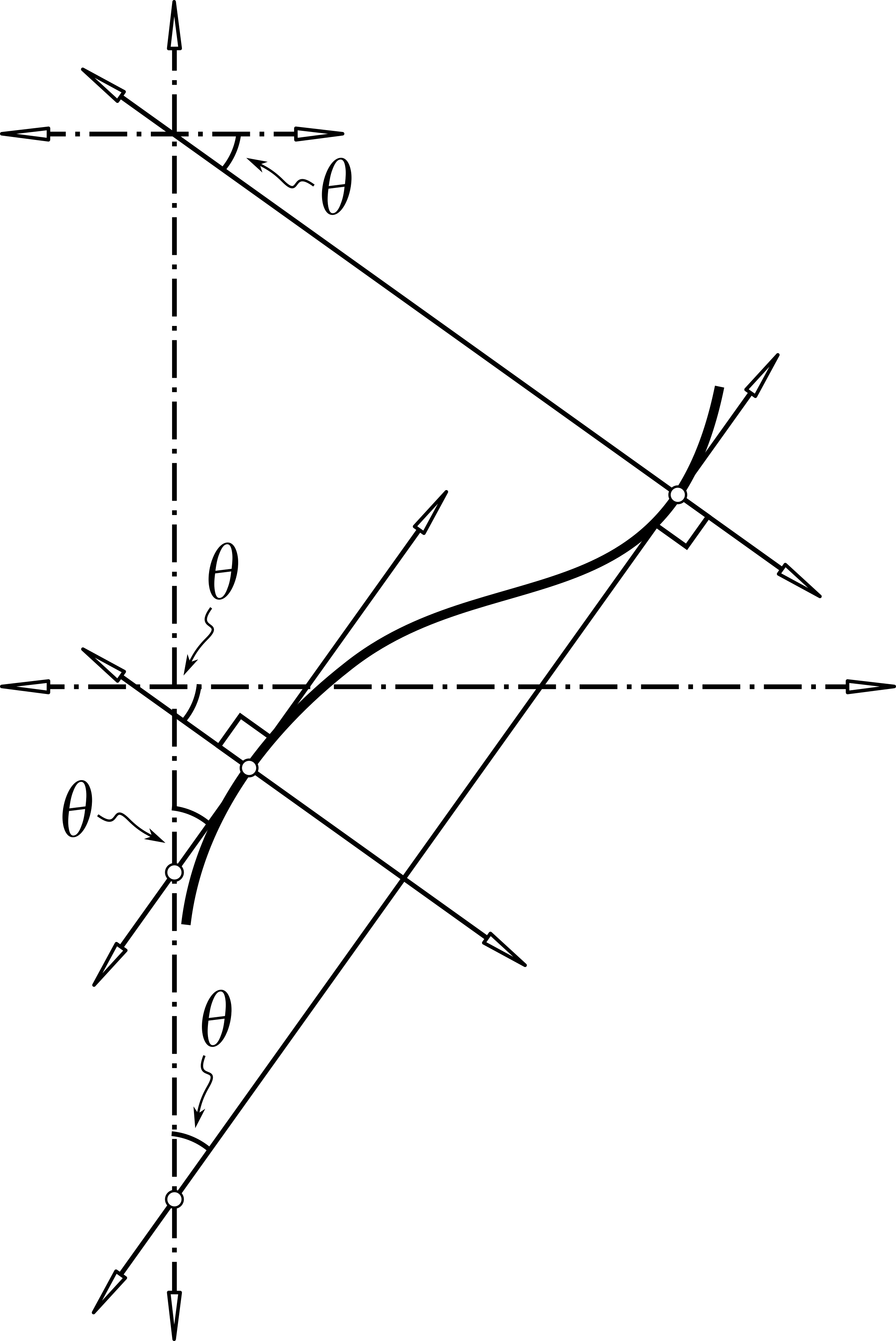}
    \end{minipage}
    \hspace{1cm}
    \begin{minipage}[m]{0.45\textwidth}
      \includegraphics[width=\textwidth]{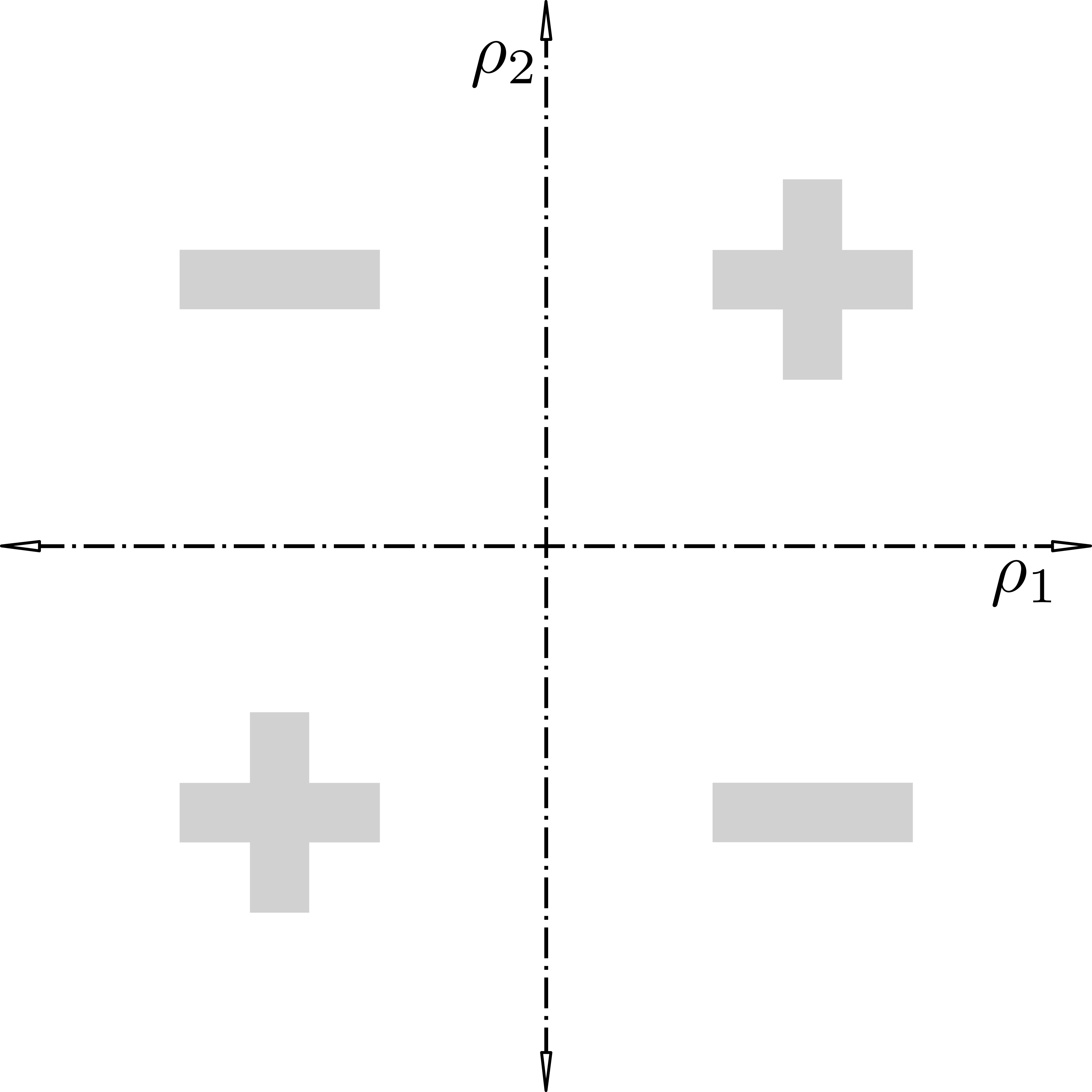}
    \end{minipage}
  }
  \caption{%
  \emph{Left}: An example of a curve for which the parametrization via its tangent lines in terms of $v \p{\theta}$ fails.  \emph{Right}: The signs of the Gau\ss ian curvatures in the quadrants of the $\rho_{1}$-$\rho_{2}$ plane.
  }
   \label{fig:v-param}
\end{figure}

The profile curve of $\Sigma$ will be parametrized via its support, as discussed in section 6.8 of \cite{Gray}; this is shown in \ref{fig:radii}, left.  In this way, the profile curve is considered as the envelope of a one-parameter family of lines, which are said to ``support'' the curve: these lines are those intersecting the axis of revolution at height $v \p{\theta}$ with an angle of $\theta \in \R$, which are tangent to the curve at $\p{r \p{\theta} , h \p{\theta}}$.
\remark{This support parametrization in terms of $v \p{\theta}$ does not preclude any allowable profile curves: the only curves for which this parametrization fails, are those where $v \p{\theta}$ fails to be a well-defined function of $\theta$; this is to say, that two different points of the profile curve have parallel tangents, corresponding to that same $\theta$, as shown in \ref{fig:v-param}, left.  This only happens when $\rho_{1}$ is infinite somewhere in between, which, by their relation \eqref{eq:aff-lin}, implies that $\rho_{2}$ must also be infinite there.  Furthermore, this can only happen if the sign of the Gau\ss ian curvature flips between those points, as $\rho_{1}$ must flip sign while $\rho_{2}$ maintains its sign. However, together, those conditions are not admitted: it is not possible for the sign of the Gau\ss ian curvature to change with $\rho_{1}$, $\rho_{2}$ going through infinity, while maintaining any such relation~\eqref{eq:aff-lin}.  To see this, note that \eqref{eq:aff-lin} is a line in the $\rho_{1}$-$\rho_{2}$ plane, which has non-zero finite slope, as $m \ne 0$ and neither $\rho_{1}$ nor $\rho_{2}$ is constant: at infinity, this line is in two non-adjacent quadrants, both of which have $\rho_{1}$, $\rho_{2}$ values that correspond to Gau\ss ian curvatures of the same sign.}

The reason for choosing this parametrization is so that it is easier to obtain expressions in terms of the principal curvature radii, which will be determined by a differential equation.  To this end, first the radii will be determined in terms of $v \p{\theta}$.  The lines forming the support, which determined by the values $\theta$ and $v \p{\theta}$, can be written in Hesse normal form as: $x, y \in \R$, such that
\[\p{\mbinom{x}{y} - \mbinom{0}{v \p{\theta}}} \cdot \mbinom{\cos \theta}{\sin \theta} = 0 .\]
Therewith, the envelope of these lines can be written as: $x, y \in \R$, such that
\[\left\{\begin{aligned}
    x \cos \theta + y \sin \theta - v \p{\theta} \sin \theta &= 0\\
    - x \sin \theta + y \cos \theta - \dot{v} \p{\theta} \sin \theta - v \p{\theta} \cos \theta &= 0
  \end{aligned}\right. \mpn{,}\]
which is satisfied by
\begin{equation}
  \left\{\begin{aligned}
    x &= - \dot{v} \p{\theta} \sin^{2} \theta \eqdef r \p{\theta}\\
    y &= v \p{\theta} + \dot{v} \p{\theta} \cos \theta \sin \theta \eqdef h \p{\theta}
  \end{aligned}\right. \mpn{.}\label{eq:env}
\end{equation}
As such, the speed of the profile curve $\p{r \p{\theta}, h \p{\theta}}$ is
\delimitershortfall=2pt
\[
  \sqrt{\p{\dot{r} \p{\theta}}^{2} + \p{\dot{h} \p{\theta}}^{2}} = \abs{\ddot{v} \p{\theta} \sin \theta + 2 \dot{v} \p{\theta} \cos \theta} \mpn{.}
\]
Computing the principal curvature radii, as per section 3C of \cite{KuehnelDiffGeom}:
\begin{align*}
  \rho_{1} \p{\theta} &= \f{\p{\p{\dot{r} \p{\theta}}^{2} + \p{\dot{h} \p{\theta}}^{2}}^{\f{3}{2}}}{\dot{r} \p{\theta} \ddot{h} \p{\theta} - \ddot{r} \p{\theta} \dot{h} \p{\theta}}\\
    &= \f{\abs{\ddot{v} \p{\theta} \sin \theta + 2 \dot{v} \p{\theta} \cos \theta}^{3}}{\p{2 \dot{v} \p{\theta} \cos \theta + \ddot{v} \p{\theta} \sin \theta}^{2}}\\
    &= \abs{\ddot{v} \p{\theta} \sin \theta + 2 \dot{v} \p{\theta} \cos \theta}\\
  \rho_{2} \p{\theta} &= \f{r \p{\theta} \sqrt{\p{\dot{r} \p{\theta}}^{2} + \p{\dot{h} \p{\theta}}^{2}}}{\dot{h} \p{\theta}}\\
    &= - \f{\dot{v} \p{\theta} \sin^{2} \theta \abs{\ddot{v} \p{\theta} \sin \theta + 2 \dot{v} \p{\theta} \cos \theta}}{2 \dot{v} \p{\theta} \cos^{2} \theta + \ddot{v} \p{\theta} \cos \theta \sin \theta}\\
    &= - \dot{v} \p{\theta} \f{\sin^{2} \theta}{\cos \theta} \sgn \p{\ddot{v} \p{\theta} \sin \theta + 2 \dot{v} \p{\theta} \cos \theta} \mpn{.}
\end{align*}
Inherent to this parametrization, there may be jumps in the normal of $\Sigma$, relative to the orientation of the profile curve in its plane.  These jumps can be seen, for example, in $\rho_{1}$ being non-negative: the normal may flip in the profile's plane so that $\rho_{1}$ maintains its sign.  In order to resolve this, the principal curvatures can be multiplied by
\[\sgn \p{- \dot{r} \p{\theta} \sin \theta} = \sgn \p{\ddot{v} \p{\theta} \sin \theta + 2 \dot{v} \p{\theta} \cos \theta} = \sgn \p{\dot{h} \p{\theta} \cos \theta} \mpn{,}\]
making
\begin{align*}
  \rho_{1} \p{\theta} &= \ddot{v} \p{\theta} \sin \theta + 2 \dot{v} \p{\theta} \cos \theta\\
  \rho_{2} \p{\theta} &= - \dot{v} \p{\theta} \f{\sin^{2} \theta}{\cos \theta} \mpn{.}
\end{align*}
Therewith, using $\rho_{2}$, it is possible also to write the coordinates of the profile curve from \eqref{eq:env} as
\[\begin{aligned}
    r \p{\theta} &= \rho_{2} \p{\theta} \cos \theta\\
    h \p{\theta} &= - \int \rho_{2} \p{\theta} \f{\cos \theta}{\sin^{2} \theta} \ \di \theta - \rho_{2} \p{\theta} \f{\cos^{2} \theta}{\sin \theta} \mpn{.}
  \end{aligned} \label{eq:r-h_rho2}\tag{$\dagger$}\]

Lastly, it remains to find the differential equations relating the principal curvature radii.  For the first of these, it is also possible to write $\rho_{1}$ in terms of $\rho_{2}$:
\allowdisplaybreaks[0]
\begin{align*}
  \rho_{2} \p{\theta} - \dot{\rho_{2}} \p{\theta} \cot \theta &= - \dot{v} \p{\theta} \f{\sin^{2} \theta}{\cos \theta} + \cot \theta \p{\ddot{v} \p{\theta} \f{\sin^{2} \theta}{\cos \theta} + \dot{v} \p{\theta} \p{2 \sin \theta + \f{\sin^{3} \theta}{\cos^{2} \theta}}\!}\\
  &= \ddot{v} \p{\theta} \sin \theta + 2 \dot{v} \p{\theta} \cos \theta\\
  &= \rho_{1} \p{\theta} \mpn{.}
\end{align*}
\allowdisplaybreaks
With their relation \eqref{eq:aff-lin}, this means that $\rho_{2}$ satisfies
\[\p{m + 1} \rho_{2} \p{\theta} - \dot{\rho_{2}} \p{\theta} \cot \theta = c \mpn{.} \label{eq:rho2_diffeq}\tag{$\ddagger$}\]
\delimitershortfall=0pt

\section{Solutions when $m = - 1$}
\label{sec:m-1}

Looking at \eqref{eq:rho2_diffeq}, it is quick to see that $m = - 1$ is a special case, as it changes the form of the differential equation.  Solving \eqref{eq:rho2_diffeq} in this case, makes
\[\rho_{2} \p{\theta} = c \log \cos \theta + J \mpn{,}\]
for some constant $J \in \R$, where $\theta \in \p{- \f{\pi}{2} , \f{\pi}{2}}$; with \eqref{eq:aff-lin}, this yields
\[\rho_{1} \p{\theta} = c \log \cos \theta + c + J \mpn{.}\]
Thus, the profile curve, for given $c, J, K \in \R$, is given by \eqref{eq:r-h_rho2} as
\begin{align*}
  r \p{\theta} &= c \cos \theta \log \cos \theta + J \cos \theta\\
  h \p{\theta} &= - c \int \f{\cos t}{\sin^{2} t} \log \cos t \;\di t - J \int \f{\cos t}{\sin^{2} t} \;\di t - \rho_{2} \p{\theta} \f{\cos^{2} \theta}{\sin \theta}\\
    &= - c \p{- \csc \theta \log \cos \theta - \int \csc t \tan t \;\di t} + J \csc \theta - \rho_{2} \p{\theta} \f{\cos^{2} \theta}{\sin \theta} + K'\\
    &= c \p{\csc \theta \log \cos \theta + \int \f{\cos t}{1 - \sin^{2} t} \;\di t} + J \csc \theta - \rho_{2} \p{\theta} \f{\cos^{2} \theta}{\sin \theta} + K'\\
    &= c \p{\csc \theta \log \cos \theta + \log \p{\sec \theta + \tan \theta}} + J \csc \theta - \rho_{2} \p{\theta} \f{\cos^{2} \theta}{\sin \theta} + K\\
    &= c \p{\sin \theta \log \cos \theta + \log \p{\sec \theta + \tan \theta}} + J \sin \theta + K \mpn{.}
\end{align*}
By construction, the normal to these curves is $n \p{\theta} = \p{- \cos \theta, - \sin \theta}$, making the evolute of that curve
\begin{align*}
  \epsilon_{- 1 , c} \p{\theta} &\defeq \p{r \p{\theta}, h \p{\theta}} + \rho_{1} \p{\theta} n \p{\theta}\\
    &= \p{-c \cos \theta, -c \sin \theta + c \log \p{\sec \theta + \tan \theta} + K} \mpn{.}
\end{align*}

The following proposition summarizes the behavior of the parameters $c , J , K \in \R$ for this case.

\begin{figure}[h!]
  \centering{%
    \begin{minipage}[m]{0.5\textwidth}
      \includegraphics[width=\textwidth]{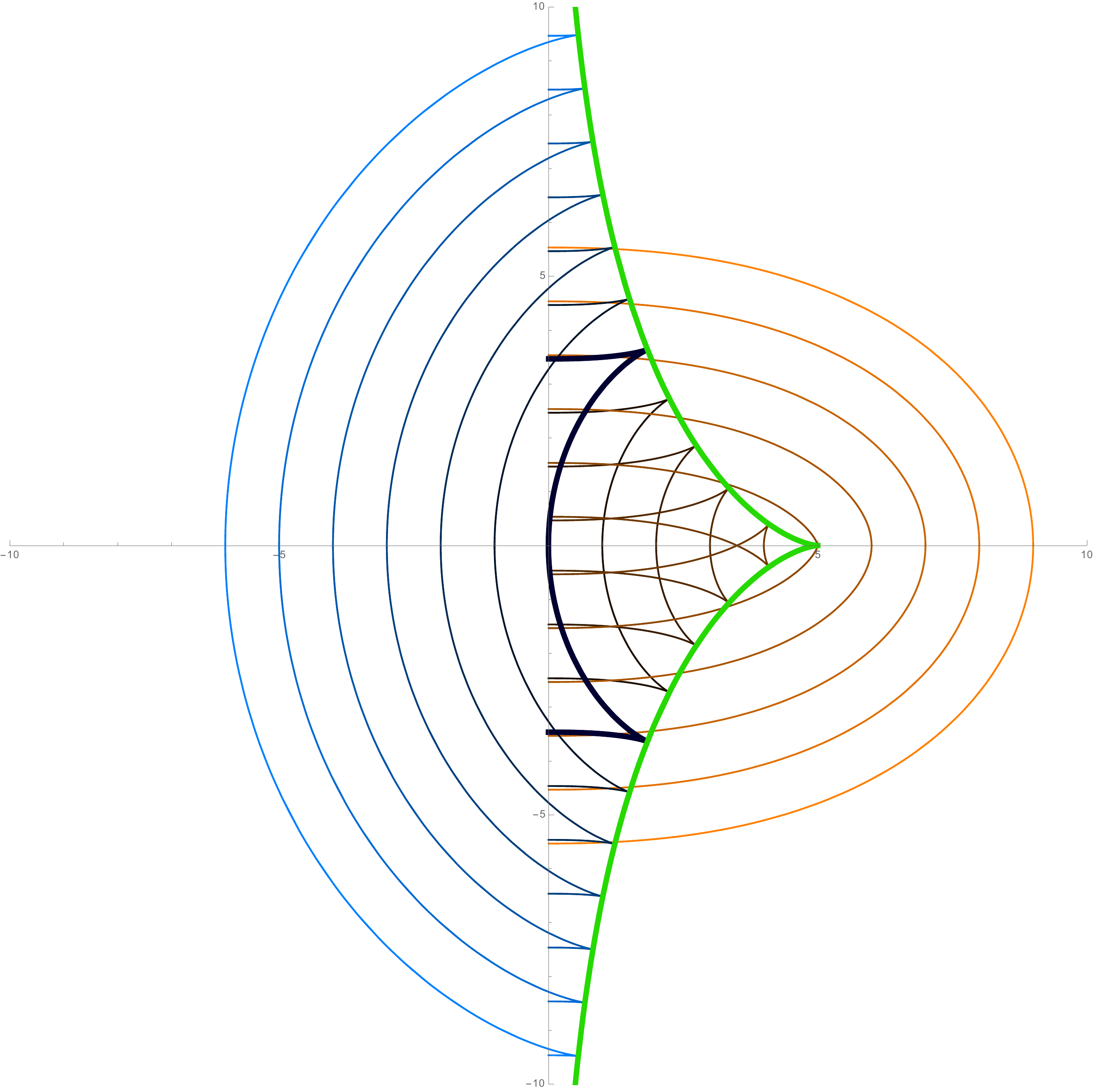}
    \end{minipage}
  }
  \caption{%
  Examples of profile curves for $m = -1$, $c = -5$, and a variety of $J$ values.  The curve with $J = 0$ is shown in bold black, while as $J$ decrease therefrom, the curves become more blue, and as it increase, the curves become more orange.  Their common evolute tractrix is shown in green.
  }
   \label{fig:examples_m-1_c-5}
\end{figure}

\prop{%
  The profile curves of the surfaces of revolution of $\p{- 1 , c}$-type for each fixed $c \in \R$ have a common evolute $\epsilon_{- 1 , c} \p{\theta}$, up to translation along the axis of revolution; in particular, when $c \ne 0$, these evolutes are tractrices.  Moreover, up to those translations, among those profile curves for $J \in \cur{0 , c}$, and among those evolutes, there is a similarity, where $c$ is the overall scaling factor.

  When $c = 0$, profile curves become circular, and the evolutes degenerate into a point along the axis of revolution.

  The integration constant $J$ parameterizes the family of normal offsets for a given $c$, see \ref{fig:examples_m-1_c-5}, while $K$ corresponds to a translations along the axis of revolution.
}

\prf{%
  Up to translations, determined by $K \in \R$, along the axis of revolution, these surfaces are those whose profile curves are involutes, determined by $J \in \R$, of the curve $\epsilon_{- 1 , c} \p{\theta}$, determined by $c \in \R$.  Moreover, as expected from the difference of the principal curvature radii being a constant, $\epsilon_{- 1 , c} \p{\theta}$ is a tractrix along the axis of revolution, when $c \ne 0$: to show this, it will be shown that the signed length $L$, along its tangent, to the axis is constant up to sign.  This $L$ is such that
  \begin{align*}
    \epsilon_{- 1 , c}^{1} \p{\theta} + L \f{\dot{\epsilon}_{-1 , c}^{1} \p{\theta}}{\abs{\dot{\epsilon}_{- 1 , c} \p{\theta}}} &= 0\\
    \xLeftrightarrow{\msp\theta \ne 0\msp} L &= \f{c \cos \theta \sqrt{\p{c \sin \theta}^{2} + \p{- c \cos \theta + c \sec \theta}^{2}}}{c \sin \theta}\\
    &= \cot \theta \sqrt{c^{2} \tan^{2} \theta}\\
    &= \abs{c} \sgn \tan \theta \mpn{,}
  \end{align*}
  and $\abs{L} = \abs{c}$ at $\theta = 0$, so $L$ has constant absolute value.  In particular, $\p{-c , 0}$ is the starting point ($\theta = 0$) of the evolute tractrix, and $\p{J , 0}$ is the starting point of involute profile curve.

  For $J \in \cur{0 , c}$, $c$ factors out of both $\p{r \p{\theta} , h \p{\theta}}$ and $\epsilon_{- 1 , c} \p{\theta}$, up to the translation by $K$, making it a scaling factor.  Thus, those curves $\p{r \p{\theta} , h \p{\theta}}$ across $c \in R$ are all similar to each other, and the $\epsilon_{- 1 , c} \p{\theta}$ are likewise.

  The curve $\epsilon_{- 1 , 0} \p{\theta}$ is degenerate, as a point along the axis of revolution, and $J$ determines the radius of the circular arc of profile curve.  This agrees with the fact that $\rho_{1} = \rho_{2}$ when $c = 0$, making the resulting surface totally umbilic.

  The fact that $J$ parameterizes the uniform normal offsets quickly follows from the normal being $\p{- \cos \theta , - \sin \theta}$ and the fact that the only terms in $r \p{\theta}$, $h \p{\theta}$ involving $J$ are respectively $J \cos \theta$, $J \sin \theta$.
}

\section{Solutions when $m \ne - 1$}
\label{sec:mne-1}

Here the general solution of \eqref{eq:rho2_diffeq} for $m \ne -1$ will be obtained.  To that end, first the homogeneous equation,
\[\p{m + 1} \rho_{2} \p{\theta} - \dot{\rho_{2}} \p{\theta} \cot \theta = c \overset{!}{=} 0 \mpn{,}\]
is solved: for some constant $J \in \R$,
\[\rho_{2} \p{\theta} = J \sec^{m + 1} \theta \mpn{.}\]
Then, therefrom a general solution can be obtained by adding the specific solution $\rho_{2} \p{\theta} = \f{c}{m + 1}$ of the inhomogeneous equation: this makes
\[\rho_{2} \p{\theta} = J \sec^{m + 1} \theta + \f{c}{m + 1}  \mpn{.}\]
As such, from \eqref{eq:aff-lin},
\[\rho_{1} \p{\theta} = - m J \sec^{m + 1} \theta + \f{c}{m + 1} \mpn{.}\]
Thus, \eqref{eq:r-h_rho2} yields the profile curve: for given $c, J, K' \in \R$,
\begin{equation}\label{eq:profile_m}
  \begin{aligned}
    r \p{\theta} &= J \sec^{m} \theta + \f{c}{m + 1} \cos \theta\\
    h \p{\theta} &= - J \int \f{\sec^{m} \theta}{\sin^{2} \theta} \ \di \theta - \f{c}{m + 1} \int \f{\cos \theta}{\sin^{2} \theta} \ \di \theta - J \f{\sec^{m - 1} \theta}{\sin \theta} - \f{c}{m + 1} \f{\cos^{2} \theta}{\sin \theta}\\
      &= - J \int \p{\f{\sec^{m} \theta}{\sin^{2} \theta} + \p{m - 1} \sec^{m} \theta - \f{\sec^{m - 2} \theta}{\sin^{2} \theta}} \ \di \theta\\
        &\qquad - \f{c}{m + 1} \int \p{\f{\cos \theta}{\sin^{2} \theta} - 2 \cos \theta - \cos \theta \cot^{2} \theta} \ \di \theta\\
      &= - m J \int \sec^{m} \theta \ \di \theta + \f{c}{m + 1} \int \cos \theta \ \di \theta\\
      &= - m J \int \sec^{m} \theta \ \di \theta + \f{c}{m + 1} \sin \theta + K' \mpn{.}
  \end{aligned}
\end{equation}
And using $\sec^{m} \theta \tan \theta = \int \sec^{m} \theta \ \di \theta + \p{m + 1} \int \sec^{m} \theta \tan^{2} \theta \ \di \theta$, its evolute is
\begin{equation}\label{evolute:m}
  \begin{aligned}
  \epsilon_{m , c} \p{\theta} &= \p{r \p{\theta}, h \p{\theta}} + \rho_{1} \p{\theta} n \p{\theta}\\
    &= \p{\p{m + 1} J \sec^{m} \theta ,\msp m \p{m + 1} J \int \sec^{m} \theta \tan^{2} \theta \ \di \theta + K'} \mpn{,}
  \end{aligned}
\end{equation}
wherein $n \p{\theta} = \p{- \cos \theta, - \sin \theta}$ is again the normal of the profile curve, by construction.  Note that, in contrast to the case when $m = -1$, this evolute is independent of $c$, and $J$ is a scaling factor, which is vice versa of the case when $m = - 1$, as shown in \ref{2.1}.  The geometric reason it is independent of $c$ follows from the fact that $c$ corresponds to the amount of the uniform normal offset.  Also, note that, when $J = 0$, the profile curve becomes a circle of radius $\f{c}{m + 1}$.

\remark{%
  When $m = 0$, it is quick to see that $\p{r \p{\theta}, h \p{\theta}}$ is a circle, as mentioned in \Aref{sec:support}, where $\p{J , K'}$ is its center, and $c$ is its radius.  Naturally, the evolute $\epsilon_{0 , c} \p{\theta}$ becomes a point, the center of that circle.
}

\section{Explicit parametrizations when $m \in \Z$}
\label{sec:alg}

By restricting $m$ to just integers values, it is possible to explicitly integrate the coordinates of the profile curves of surfaces of type $\p{m , c}$-type from \eqref{eq:profile_m}.  With the explicit representations of their profile curves, it is possible then also to determine their algebraicity.  It is well-known that a surface of revolution with an algebraic profile curve is itself algebraic, and that its algebraic degree is at most twice the degree of its profile curve.

Starting off, the following proposition proves the conjecture from \cite{Pottmann} that surfaces of $\p{m , 0}$-type are algebraic for even $m \in \Z_{\ge 0}$, by showing it for $\p{m , c}$-type.

\prop{%
  For even $m \in \Z_{\ge 0}$, the profile curve $\p{r \p{\theta}, h \p{\theta}}$ from \eqref{eq:profile_m} is algebraic for any $c , J , K \in \R$.  Moreover, when $c = 0$, it is an algebraic curve of degree $m$, while when $c \ne 0$, its degree is $2 \p{m + 1}$.
}

\prf{%
  As discussed in \ref{3.1}, when $m = 0$, the profile curve is simply a circle, which is, of course, algebraic.

  For even values $m > 0$,
  \begin{align*}
    r \p{\theta} &= J \sec^{m} \theta + \f{c}{m + 1} \cos \theta\\
    h \p{\theta} &= - m J \int \sec^{m} \theta \ \di \theta + \f{c}{m + 1} \sin \theta + K'\\
      &= - m J \int \p{1 + u^{2}}^{\f{m - 2}{2}} \ \di u + \f{c}{m + 1} \sin \theta + K' &&\mathrlap{\hspace{-2.5ex}\mnote{10ex}{$u = \tan \theta$}}\\
      &= - m J \sum_{j = 0}^{\f{m - 2}{2}} \mbinom{\f{m - 2}{2}}{j} \f{\tan^{2 j + 1} \theta}{2 j + 1} + \f{c}{m + 1} \sin \theta + K \mpn{.}
  \end{align*}
  As $\p{r \p{\theta + \pi}, h \p{\theta + \pi}} = \p{r \p{\theta}, - h \p{\theta}}$, and there is an overall period of $2 \pi$, it suffices to show that the restriction of the profile curve to $\theta \in \ioo{- \f{\pi}{2}, \f{\pi}{2}}$ is algebraic, in order to show that the entire profile curve is algebraic.  To this end, reparametrizing with $\theta \p{t} = \arctan t$ for $t \in \R$,
  \delimitershortfall=5pt
  \begin{align*}
    r \circ \theta \p{t} &= J \p{t^{2} + 1}^{\f{m}{2}} + \f{c}{m + 1} \f{1}{\sqrt{t^{2} + 1}}\\
    h \circ \theta \p{t} &= - m J \sum_{j = 0}^{\f{m - 2}{2}} \mbinom{\f{m - 2}{2}}{j} \f{t^{2 j + 1}}{2 j + 1} + \f{c}{m + 1} \f{t}{\sqrt{t^{2} + 1}} + K \mpn{.}
  \end{align*}
  Wherefrom it is possible to obtain polynomial representations in, respectively, $\R \s{r , t}$ and $\R \s{h , t}$:
  \begin{align*}
    \p{t^{2} + 1} \p{r - J \p{t^{2} + 1}^{\f{m}{2}}}^{2} - \p{\f{c}{m + 1}}^{2} &= 0\\
    \p{t^{2} + 1} \p{h + m J \sum_{j = 0}^{\f{m - 2}{2}} \mbinom{\f{m - 2}{2}}{j} \f{t^{2 j + 1}}{2 j + 1} - K}^{2} - \p{\f{c}{m + 1}}^{2} t^{2} &= 0 \mpn{,}
  \end{align*}
  showing the profile curve is algebraic.

  As the restriction of the profile curve to $\theta \in \p{- \f{\pi}{2} , \f{\pi}{2}}$ is injective, finding its degree is the same as finding the number of points, at which the curve intersects a generic line $A r + B h = C$ for $A , B , C \in \C$.  When $c \ne 0$, this intersection is given by solutions to the following polynomial
  \begin{align*}
    &\p{t^{2} + 1} \p{A J \p{t^{2} + 1}^{\f{m}{2}} + B \p{- m J \sum_{j = 0}^{\f{m - 2}{2}} \mbinom{\f{m - 2}{2}}{j} \f{t^{2 j + 1}}{2 j + 1} + K} - C}^{2}\\
    &\hspace{10ex}- \p{\f{\p{A + B t} c}{m + 1}}^{2} = 0 \mpn{,}
  \end{align*}
  which has degree $2 \p{m + 1}$ in $t$.  And when $c = 0$, this is simply
  \[A J \p{t^{2} + 1}^{\f{m}{2}} + B \p{- m J \sum_{j = 0}^{\f{m - 2}{2}} \mbinom{\f{m - 2}{2}}{j} \f{t^{2 j + 1}}{2 j + 1} + K} = C \mpn{,}\]
  which is of degree $m$ in $t$.
}

\prop{%
  For $m = - 1$, the profile curve $\p{r \p{\theta}, h \p{\theta}}$ from \eqref{eq:profile_m} is algebraic for $c = 0$ and any $J , K \in \R$, while for odd $m \in \Z_{< -1}$, it is algebraic for any $c , J , K \in \R$.  Moreover, for any odd $m \in \Z_{< 0}$, the algebraic degree of this curve is $- 2 m$.
}

\prf{%
  As discussed in \Aref{sec:m-1} for the case $m = - 1$, the profile curve is a circular arc, when $c = 0$: as such, it is algebraic for any $J , K \in \R$, and of degree $2$.

  Letting $\hat{m} \defeq - m$ as to simplify notation: for negative, odd values of $m < - 1$,
  \delimitershortfall=5pt
  \begin{align*}
    r \p{\theta} &= J \cos^{\hat{m}} \theta + \f{c}{1 - \hat{m}} \cos \theta\\
    h \p{\theta} &= \hat{m} J \int \cos^{\hat{m}} \theta \ \di \theta + \f{c}{1 - \hat{m}} \sin \theta + K'\\
      &= \hat{m} J \int \p{1 - u^{2}}^{\f{\hat{m} - 1}{2}} \ \di u + \f{c}{1 - \hat{m}} \sin \theta + K' &&\mathrlap{\hspace{-2.5ex}\mnote{10ex}{$u = \sin \theta$}}\\
      &= \hat{m} J \sum_{j = 0}^{\f{\hat{m} - 1}{2}} \p{-1}^{j} \mbinom{\f{\hat{m} - 1}{2}}{j} \f{\sin^{2 j + 1} \theta}{2 j + 1} + \f{c}{1 - \hat{m}} \sin \theta + K \mpn{.}
  \end{align*}
  Because of its periodicity, it suffices to show that, for $\theta \in \p{- \pi , \pi}$, the profile curve is algebraic.  To that end, reparametrizing with $\theta = \arcsin t$ yields
  \begin{align*}
    r \p{\theta} &= J \p{1 - t^{2}}^{\f{\hat{m}}{2}} + \f{c}{1 - \hat{m}} \sqrt{1 - t^{2}}\\
    h \p{\theta} &= \hat{m} J \sum_{j = 0}^{\f{\hat{m} - 1}{2}} \p{-1}^{j} \mbinom{\f{\hat{m} - 1}{2}}{j} \f{t^{2 j + 1}}{2 j + 1} + \f{c}{1 - \hat{m}} t + K \mpn{.}
  \end{align*}
  Wherefrom polynomial representations in, respectively, $\R \s{r , t}$ and $\R \s{h , t}$, can be found:
  \begin{align*}
    r^{2} - \p{1 - t^{2}} \p{J \p{1 - t^{2}}^{\f{\hat{m} - 1}{2}} + \f{c}{1 - \hat{m}}}^{2} &= 0\\
    h - \hat{m} J \sum_{j = 0}^{\f{\hat{m} - 1}{2}} \p{-1}^{j} \mbinom{\f{\hat{m} - 1}{2}}{j} \f{t^{2 j + 1}}{2 j + 1} - \f{c}{1 - \hat{m}} t - K &= 0 \mpn{,}
  \end{align*}
  showing that the profile curve is algebraic for any $c , J , K \in \R$.

  As the restriction of the profile curve to $\theta \in \p{- \pi , \pi}$ is injective, finding its degree is the same as finding the number of points, at which the curve intersects a generic line $A r + B h = C$ for $A , B , C \in \C$.  When $c \ne 0$, this intersection is given by solutions to the following polynomial
  \begin{multline*}
    \p{B \p{\hat{m} J \sum_{j = 0}^{\f{\hat{m} - 1}{2}} \p{-1}^{j} \mbinom{\f{\hat{m} - 1}{2}}{j} \f{t^{2 j + 1}}{2 j + 1} + \f{c}{1 - \hat{m}} t + K} - C}^{2}\\
    - \p{A \p{J \p{1 - t^{2}}^{\f{\hat{m} - 1}{2}} + \f{c}{1 - \hat{m}}}}^{2} \p{1 - t^{2}} = 0 \mpn{,}\!\phantom{,}
  \end{multline*}
  which is of degree $2 \hat{m}$ in $t$.  And when $c = 0$, this is simply
  \begin{multline*}
    \p{B \p{\hat{m} J \sum_{j = 0}^{\f{\hat{m} - 1}{2}} \p{-1}^{j} \mbinom{\f{\hat{m} - 1}{2}}{j} \f{t^{2 j + 1}}{2 j + 1} + K} - C}^{2}\\
    - \p{A J \p{1 - t^{2}}^{\f{\hat{m}}{2}}}^{2} = 0 \mpn{,}\!\phantom{,}
  \end{multline*}
  which is also of degree $2 \hat{m}$ in $t$.
}

\prop{%
  For odd $m \in \Z_{> 0}$, the profile curve $\p{r \p{\theta}, h \p{\theta}}$ from \eqref{eq:profile_m} is transcendental for any $c , K \in \R$, $J \in \R^{\times}$.
}

\prf{%
  For $m = 1$,
  \begin{align*}
    r \p{\theta} &= J \sec \theta + \f{c}{2} \cos \theta\\
    h \p{\theta} &= - J \int \sec \theta \ \di \theta + \f{c}{2} \sin \theta + K'\\
      &= - J \log \abs{\sec \theta + \tan \theta} + \f{c}{2} \sin \theta + K \mpn{.}
  \end{align*}
  And, for odd values $m > 1$,
  \delimitershortfall=5pt
  \begin{align*}
    r \p{\theta} &= J \sec^{m} \theta + \f{c}{m + 1} \cos \theta\\
    h \p{\theta} &= - m J \int \sec^{m} \theta \ \di \theta + \f{c}{m + 1} \sin \theta + K'\\
      &= - m J \p{\f{\sec^{m - 2} \theta \tan \theta}{m - 1} + \f{m - 2}{m - 1} \int \sec^{m - 2} \theta \ \di \theta} + \f{c}{m + 1} \sin \theta + K''\\
      &= - m J \p{\sum_{j = 1}^{\f{m - 1}{2}} \f{\p{m - 2}!!}{\p{m - 2j}!!} \f{\p{m - 2 j + 1}!!}{\p{m - 1}!!} \f{\sec^{m - 2 j} \theta \tan \theta}{m - 2 j + 1} \right.\\
      &\hspace{15ex}\left. \vphantom{\sum_{j = 1}^{\f{m - 1}{2}}} + \f{1}{2} \log \abs{\sec \theta + \tan \theta}} + \f{c}{m + 1} \sin \theta + K \mpn{.}
  \end{align*}
  \delimitershortfall=0pt

  As discussed in \Aref{sec:m-1}, $c \in \R$ corresponds to a uniform normal offset of the profile curve.  Therewith, it suffices to show that the profile curves for $c = 0$ are transcendental, in order to show that they are transcendental for any $c \in \R$, as uniform normal offsets preserve algebraicity.  To this end, restrict to $\theta \in \ioo{0 , \f{\pi}{2}}$, and consider the reparametrizations of $\theta = \arcsec t$ for $t \in \R_{> 1}$: for $m = 1$,
  \begin{align*}
    r \p{\theta} &= J t\\
    h \p{\theta} &= - J \log \p{t + \sqrt{t^{2} - 1}} + K \mpn{;}
  \end{align*}
  for odd $m > 1$,
  \begin{align*}
    r \p{\theta} &= J t^{m}\\
    h \p{\theta} &= - m J \p{\sum_{j = 1}^{\f{m - 1}{2}} \f{\p{m - 2}!!}{\p{m - 2j}!!} \f{\p{m - 2 j + 1}!!}{\p{m - 1}!!} \f{t^{m - 2 j} \sqrt{t^{2} - 1}}{m - 2 j + 1} \right.\\
      &\hspace{15ex}\left. \vphantom{\sum_{j = 1}^{\f{m - 1}{2}}} + \f{1}{2} \log \p{t + \sqrt{t^{2} - 1}}} + K \mpn{.}
  \end{align*}
  Then, the profile curves are of the following form: for some $\psi_{m} \p{t} \in \Q \s{t}$,
  \begin{multline*}
    \p{R \p{t} , H \p{t}} \defeq\\
    \p{J t^{m} \msp,\msp \subnote{- m J \sqrt{t^{2} - 1} \psi_{m} \p{t}}{20ex}{$\eqdef \sqrt{t^{2} - 1} \ A \p{t}$} \;-\; \subnote{\f{m J}{2} \log \p{t + \sqrt{t^{2} - 1}}}{15ex}{$\eqdef - B \p{t}$}\,}\ubfixvsp \mpn{,}\!\phantom{,}
  \end{multline*}
  wherein $A \p{t} \in \R \s{t}$.  Assuming that this curve is algebraic, there is a $P \p{x , y} \in \C \s{x , y}$ such that
  \[P \p{R \p{t} , H \p{t}} = 0 \mpn{.}\]
  In particular, this means that there are $c_{j k} \in \C$ so that
  \[\sum_{j , k} c_{j k} \ t^{m j} \p{\sqrt{t^{2} - 1} \ A \p{t} + B \p{t}}^{k} = 0 \mpn{,}\]
  which is of the form
  \begin{align*}
    \sqrt{t^{2} - 1} \ F \p{t , B \p{t}} + G \p{t , B \p{t}} &= 0\\
    \implies \p{t^{2} - 1} \p{F \p{t , B \p{t}}}^{2} - \p{G \p{t , B \p{t}}}^{2} &= 0 \mpn{,}
  \end{align*}
  wherein there are some $F \p{x , y} , G \p{x , y} \in \C \s{x , y}$.  Now, consider the limit as $t \to \infty$ of that last expression: looking at the dominating terms,
  \begin{multline*}
    \p{t^{2} - 1} \p{F \p{t , B \p{t}}}^{2} - \p{G \p{t , B \p{t}}}^{2} \ \sim\\
    D_{0} t^{N} + D_{1} t^{N'} \p{\log \circ f \p{t}}^{M'} + D_{2} \p{\log \circ f \p{t}}^{M} \mpn{,}\!\phantom{,}
  \end{multline*}
  for some  $D_{j} \in \C$ not all zero, some $N, N', M, M' \in \Z_{\ge 0}$ not all zero, and with $f \p{t} \defeq t + \sqrt{t^{2} - 1}$.  As such, it suffices to show that this goes to infinity in this limit.  To this end, note that
  \begin{align*}
      t^{U} &= o \p{t^{U} \p{\log \circ f \p{t}}^{V'}}\\
      t^{U} \p{\log \circ f \p{t}}^{V} &= o \p{t^{U + 1}} \mpn{,}
  \end{align*}
  for any $U , V \in \Z_{\ge 0}$ and any $V' \in \Z_{> 0}$.  The former is quickly obtained as the logarithm is unbounded, and the latter is derived as follows: By definition, this means that, for any $\epsilon > 0$, there exists an $L \ge 0$ large-enough so that, for any $t \ge L$,
  \begin{align*}
    t^{U} \p{\log \p{\sqrt{t^{2} - 1} + q}}^{V} &\overset{!}{\le} \epsilon t^{U + 1}\\
    \iff \sqrt{t^{2} - 1} + q &\overset{!}{\le} \exp \p{\epsilon t^{\f{1}{V}}}\\
    \iff -1 &\overset{!}{\le} \exp \p{2 \epsilon t^{\f{1}{V}}} - 2 q \exp \p{\epsilon t^{\f{1}{V}}}\\
    \iff 0 &\overset{!}{\le} \f{\exp \p{2 \epsilon t^{\f{1}{V}}}}{2 q \exp \p{\epsilon t^{\f{1}{V}}} - 1} \mpn{.}
  \end{align*}
  Finding such an $L$ is equivalent to showing that the right-hand side goes to $\infty$, as $q \to \infty$.  Proceeding thusly, applying L'H\^{o}pital's rule successively, to handle the indeterminate form $\f{\infty}{\infty}$,
  \begin{align*}
    \lim_{q \to \infty} \f{\exp \p{2 \epsilon t^{\f{1}{V}}}}{2 q \exp \p{\epsilon t^{\f{1}{V}}} - 1} &= \lim_{q \to \infty} \f{2 \epsilon \p{t^{\f{1}{V}}}' \exp \p{2 \epsilon t^{\f{1}{V}}}}{2 \exp \p{\epsilon t^{\f{1}{V}}} + 2 t \epsilon \p{t^{\f{1}{V}}}' \exp \p{\epsilon t^{\f{1}{V}}}}\\
      &\hspace{-4ex}= \lim_{q \to \infty} \f{\epsilon \p{t^{\f{1}{V}}}' \exp \p{\epsilon t^{\f{1}{V}}}}{1 + t \epsilon \p{t^{\f{1}{V}}}'}\\
      &\hspace{-4ex}= \lim_{q \to \infty} \f{\exp \p{\epsilon t^{\f{1}{V}}} \p{\epsilon \p{t^{\f{1}{V}}}'' + \p{\epsilon \p{t^{\f{1}{V}}}'}^{2}}}{\epsilon \p{t^{\f{1}{V}}}' + t \epsilon \p{t^{\f{1}{V}}}''}\\
      &\hspace{-4ex}= \lim_{q \to \infty} \f{\exp \p{\epsilon t^{\f{1}{V}}} \p{1 - V + \epsilon t^{\f{1}{V}}}}{t}\\
      &\hspace{-4ex}= \lim_{q \to \infty} \exp \p{\epsilon t^{\f{1}{V}}} \p{\epsilon \p{t^{\f{1}{V}}}' \p{1 - V + \epsilon t^{\f{1}{V}}} + \epsilon \p{t^{\f{1}{V}}}'}\\
      &\hspace{-4ex}= \infty \mpn{.}
  \end{align*}
  Therewith, as was to be shown,
  \[t^{U} \p{\log \circ f \p{t}}^{V} = o \p{t^{U + 1}} \mpn{,}\]
  and in particular,
  \[\p{\log \circ f \p{t}}^{V} = \mathrm{o} \p{t}\]
  for any $V \in \Z_{\ge 0}$.  Therefore, as $t \to \infty$, it must be that
  \[\abs{D_{0} t^{N} + D_{1} t^{N'} \p{\log \circ f \p{t}}^{M'} + D_{2} \p{\log \circ f \p{t}}^{M}} \to \infty\]
  regardless of the $D_{j} \in \C$ not all zero, and the $N, N', M, M' \in \Z_{\ge 0}$ not all zero.

  Therefrom, a contradiction is found to the existence of such a $P \p{x , y} \in \C \s{x , y}$, making the profile curves $\p{R \p{t}, H \p{t}}$ transcendental for $c = 0$ and any $K \in \R$, $J \in \R^{\times}$, and thusly, any $c , K \in \R$, $J \in \R^{\times}$.
}

\prop{%
  For even $m \in \Z_{< 0}$, the profile curve $\p{r \p{\theta}, h \p{\theta}}$ from \eqref{eq:profile_m} is transcendental for any $c , K \in \R$, $J \in \R^{\times}$
}

\prf{%
  Let $\hat{m} \defeq - m$ as to simplify notation.  Then, using the following identities,
  \[
    \cos^{\hat{m}} \theta = \f{1}{2^{\hat{m}}} \mbinom{\hat{m}}{\f{\hat{m}}{2}} + \f{1}{2^{\hat{m} - 1}} \sum_{j = 0}^{\f{\hat{m} - 2}{2}} \mbinom{\hat{m}}{j} \cos \p{\hat{m} - 2 j} \theta
  \]
  and
  \[
    \sin M \theta = \sum_{\substack{k = 0\\k \text{ odd}}}^{M} \p{-1}^{\f{k - 1}{2}} \mbinom{M}{k} \cos^{M - k} \theta \sin^{k} \theta \mpn{,}
  \]
  it is possible to write the profile curves as
  \begin{align*}
    r \p{\theta} &= J \cos^{\hat{m}} \theta + \f{c}{1 - \hat{m}} \cos \theta\\
    h \p{\theta} &= \hat{m} J \int \cos^{\hat{m}} \theta \ \di \theta + \f{c}{1 - \hat{m}} \sin \theta + K'\\
      &\hspace{-6ex}= \hat{m} J \p{\f{\theta}{2^{\hat{m}}} \mbinom{\hat{m}}{\f{\hat{m}}{2}} + \f{1}{2^{\hat{m} - 1}} \sum_{j = 0}^{\f{\hat{m} - 2}{2}} \mbinom{\hat{m}}{j} \f{\sin \p{\hat{m} - 2 j} \theta}{\hat{m} - 2 j}} + \f{c}{1 - \hat{m}} \sin \theta + K\\
      &\hspace{-6ex}= \hat{m} J \p{\f{\theta}{2^{\hat{m}}} \mbinom{\hat{m}}{\f{\hat{m}}{2}} + \f{1}{2^{\hat{m} - 1}} \sum_{j = 0}^{\f{\hat{m} - 2}{2}} \mbinom{\hat{m}}{j} \sum_{\substack{k = 0\\k \text{ odd}}}^{\hat{m} - 2 j} \p{-1}^{\f{k - 1}{2}} \mbinom{\hat{m} - 2 j}{k} \f{\cos^{\hat{m} - 2 j - k} \theta \sin^{k} \theta}{\hat{m}  - 2 j}}\\
        &\hspace{9ex}+ \f{c}{1 - \hat{m}} \sin \theta + K \mpn{.}
  \end{align*}

  As discussed in \Aref{sec:m-1}, $c \in \R$ corresponds to a uniform normal offset of the profile curve.  Therewith, it suffices to show that the profile curves for $c = 0$ are transcendental, in order to show that they are transcendental for any $c \in \R$, as uniform normal offsets preserve algebraicity.  To this end, restrict to $\theta \in \ioo{0 , \pi}$, and consider the reparametrizations of $\theta = \arccos t$ for $t \in \ioo{-1 , 1}$:
  \allowdisplaybreaks[0]
  \begin{align*}
    r \p{\theta} &= J t^{\hat{m}}\\
    h \p{\theta} &= \hat{m} J \p{\f{\arccos t}{2^{\hat{m}}} \mbinom{\hat{m}}{\f{\hat{m}}{2}} \vphantom{\sum_{\substack{k = 0\\k \text{ odd}}}^{\hat{m} - 2 j}} \right.\\
    &\hspace{-1ex}\left. + \f{1}{2^{\hat{m} - 1}} \sum_{j = 0}^{\f{\hat{m} - 2}{2}} \mbinom{\hat{m}}{j} \sum_{\substack{k = 0\\k \text{ odd}}}^{\hat{m} - 2 j} \p{-1}^{\f{k - 1}{2}} \mbinom{\hat{m} - 2 j}{k} \f{t^{\hat{m} - 2 j - k} \p{1 - t^{2}}^{\f{k}{2}}}{\hat{m}  - 2 j}} + K \mpn{.}
  \end{align*}
  \allowdisplaybreaks
  As $\arccos{t}$ is not algebraic, it follows that there cannot exist a $P \p{x , y} \in \C \s{x , y}$ such that $P \p{r \p{\theta} , h \p{\theta}} = 0$.  Therefore, the profile curves are not algebraic for $c = 0$ and $K \in R$, $J \in \R^{\times}$, and thusly, for any $c , K \in \R$ and $J \in \R^{\times}$.
}

\section{Constant-angle parametrizations when $c = 0$}
\label{sec:param}

In this section, the special parametrizations of surfaces of $\p{m , 0}$-type will be found.  Surfaces of this type are those that, via \eqref{eq:aff-lin}, have a constant ratio of principal curvatures,
\[\f{\rho_{1}}{\rho_{2}} = - m \mpn{.}\]
Moreover, when $m > 0$, this means that, at every point, the asymptotic lines meet at a constant angle determined by $m$; this is discussed in section 4 of \cite{Jimenez-Mueller-Pottmann}, for example.  As such, first the asymptotic curves of surfaces with $m > 0$ will be determined, wherefrom their asymptotic parametrization can be formed: this will be done by finding a differential equation for the asymptotic curves in terms of their Frenet-Serret frame, an idea from the manuscript \cite{Pottmann}.  Lastly, as the argument does not hinge on the sign of $m$, it will be seen that a constant-angle parametrization for surfaces with $m < 0$ is therewith simultaneously found.  The results about asymptotic curves in this section up to \ref{5.3} elaborate on the results that culminate in Theorem 1 of \cite{Pottmann}.

Let $\alpha \colon \R \fdtr \Sigma$ be an asymptotic curve on $\Sigma \subset \R^{3}$, which will be assumed not to be a geodesic in $\R^{3}$, so that its Frenet-Serret frame exists almost everywhere; then, it readily follows from its definition, that the osculating plane at $\alpha \p{t}$ exists almost everywhere, and agrees with the tangent plane of $\Sigma$ there.  Outside of possible inflection points, it is possible to generate the Frenet-Serret frame $\cur{\ee{1}, \ee{2}, \ee{3}}$ along $\alpha$, where $\ee{1}$ is the unit tangent of $\alpha$, $\ee{2}$ is the unit normal of $\alpha$, and $\ee{3} \defeq \ee{1} \times \ee{2}$ is the binormal of $\alpha$.  Moreover, as the osculating plane agrees with the tangent plane, it follows that $\ee{1} \p{t}$ and $\ee{2} \p{t}$ span the tangent plane of $\Sigma$ at $\alpha \p{t}$, and $\ee{3} \p{t}$ is the surface normal at $\alpha \p{t}$.  Since $\Sigma$ is assumed to have a constant ratio of principal curvatures, it follows that the asymptotic directions meet the parallel circles of $\Sigma$ at a constant angle; thus, it is possible to express the (unit) tangent $\ff{1} \p{t}$ of the parallel circle going through $\alpha \p{t}$ as
\begin{equation}
  \ff{1} \p{t} \defeq \cos \tau \ee{1} \p{t} + \sin \tau \ee{2} \p{t} \mpn{,}\label{eq:ff1}
\end{equation}
wherein $\tau \in \ioo{0, \f{\pi}{2}}$ is a constant determined by $m = \tan^{2} \tau$.

As $\Sigma$ is assumed to be immersed in $\R^{3}$, it is possible to consider, in a canonical way, spherical images of the $\ee{j}$ and $\ff{1}$, and without ambiguity, consider them as curves $\R \fdtr \Sph^{2}$: they will all be considered on $\Sph^{2}$ relative to the axis of rotation $\ell$ of $\Sigma$.  It will be shown in the following proposition that these spherical curves are of particular types.

\prop{%
  On surfaces of $\p{m , 0}$-type for $m > 0$, the curves $\ee{1}, \ee{3} \colon \R \fdtr \Sph^{2}$, as defined just above, are respectively a spherical tractrix and a loxodrome.
}

\begin{figure}[h!]
  \centering{%
    \begin{minipage}[m]{0.3\textwidth}
      \includegraphics[height=\textwidth]{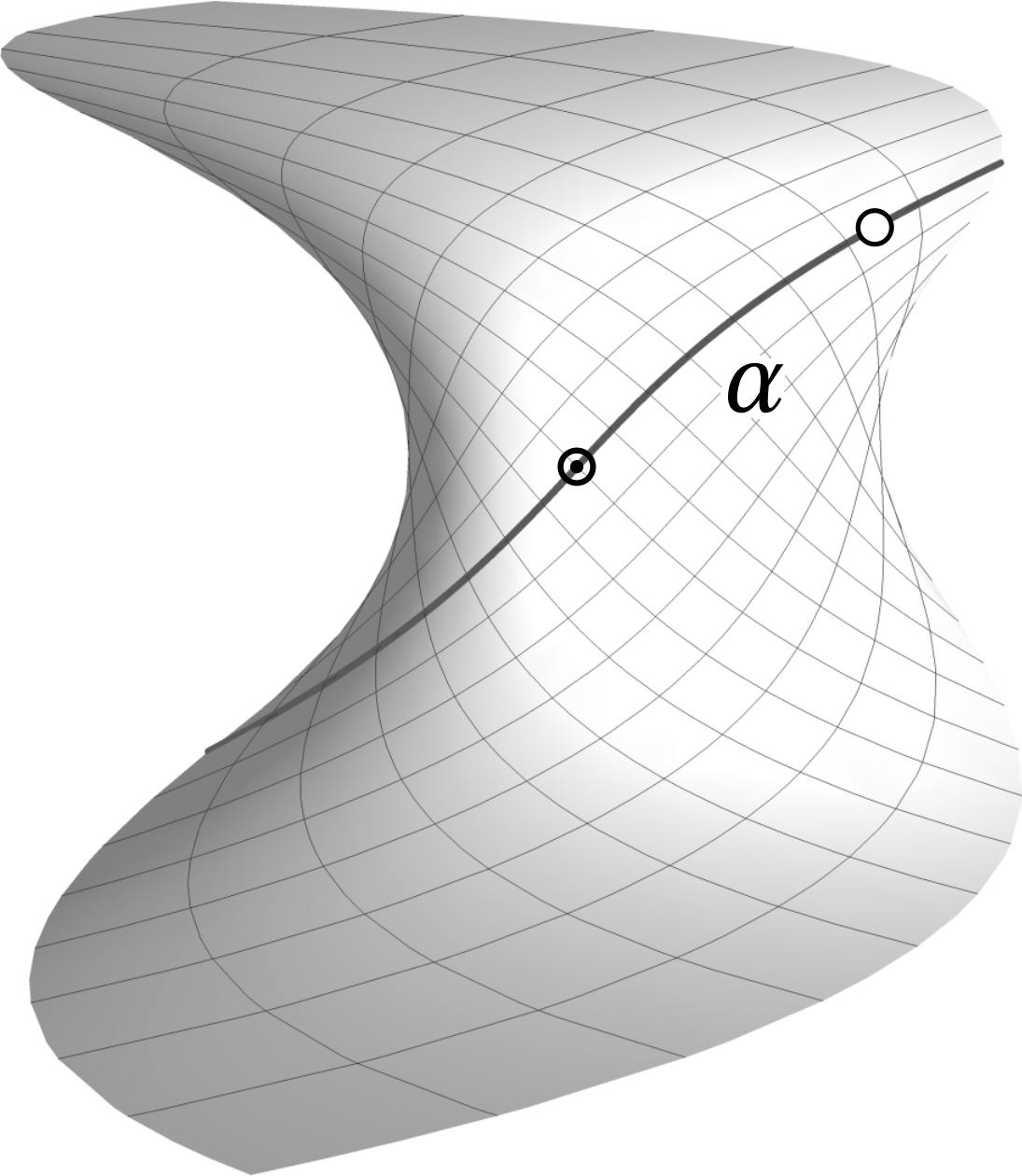}
    \end{minipage}
    \hspace{-0.05\textwidth}
    \hfill
    \begin{minipage}[m]{0.3\textwidth}
      \includegraphics[width=\textwidth]{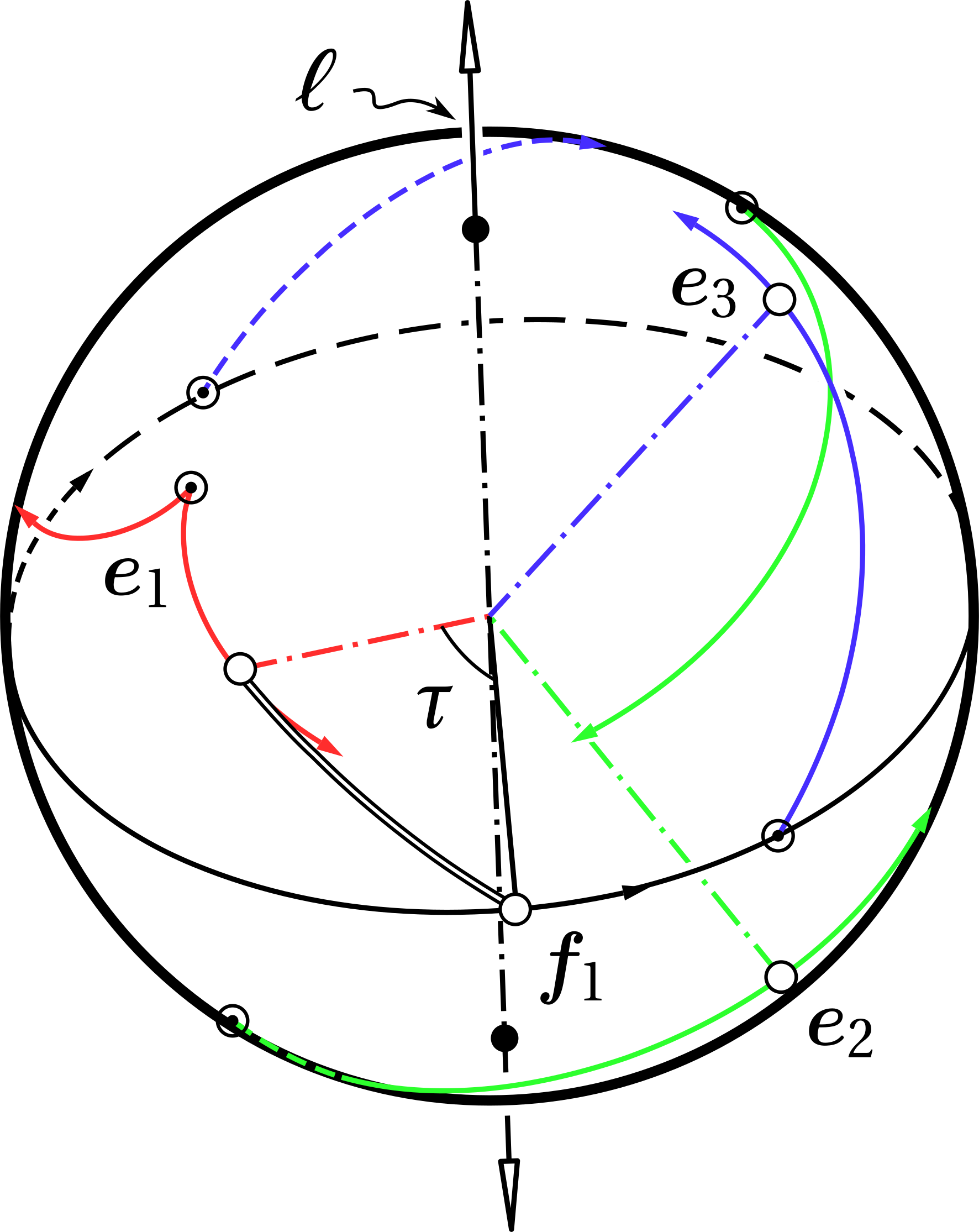}
    \end{minipage}
    \hfill
    \raisebox{0.0125\textwidth}{%
      \begin{minipage}[m]{0.325\textwidth}
        \includegraphics[width=\textwidth]{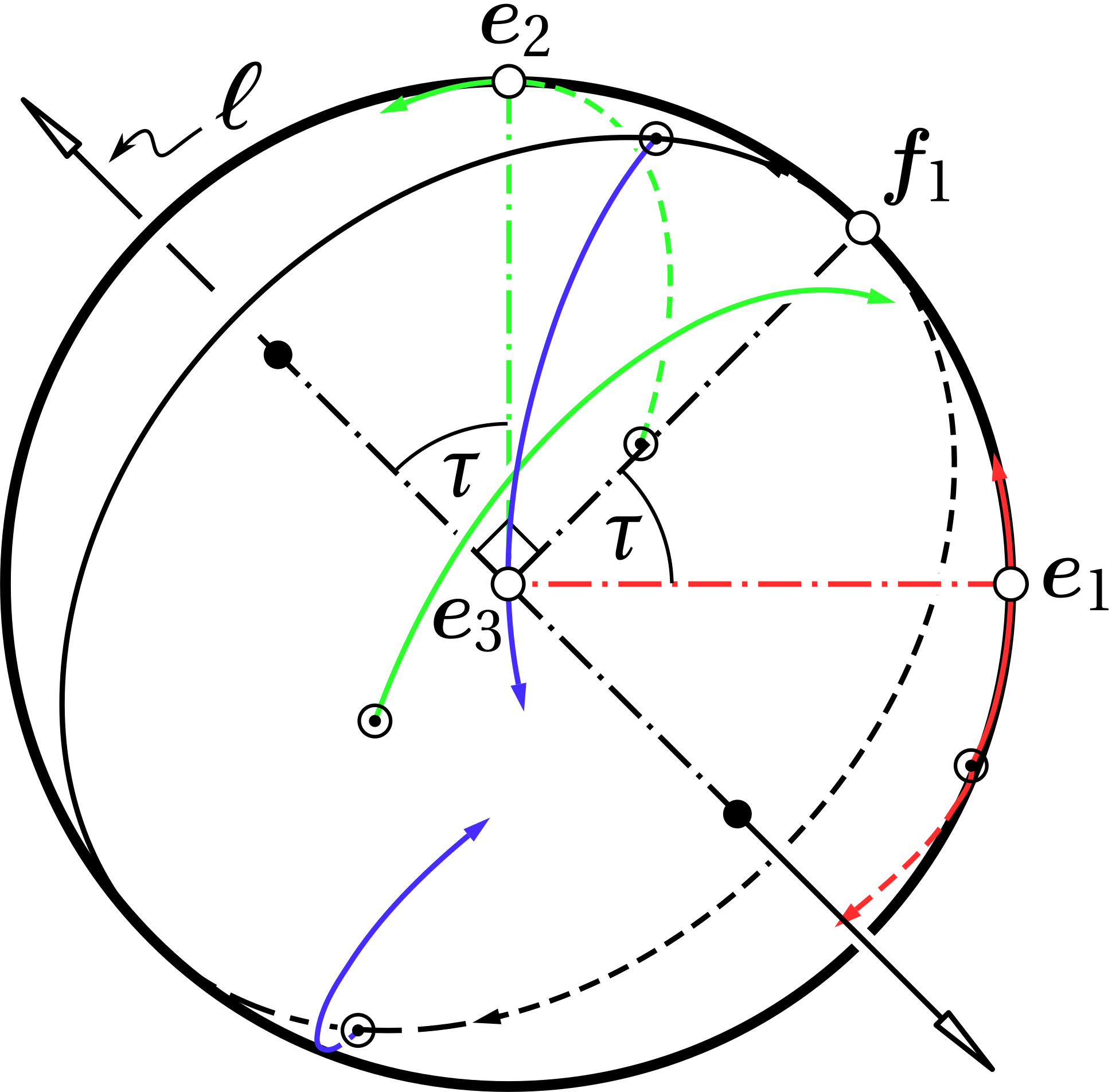}
      \end{minipage}%
    }%
  }
  \caption{%
    \emph{Left}: Catenoid with asymptotic curve along it, where the point marked with a bullseye is at the neck of $\Sigma$, and the other marked point corresponds to the frame shown in the following images of this figure.
    \emph{Center}: The spherical image of the frame along that curve, along with the curve of $\ff{1}$.  Note that the cusp in the curve $\ee{1}$, and correspondingly, the discontinuity in each of the curves $\ee{2}$ and $\ee{3}$ occur at the neck of the catenoid.
    \emph{Right}: The projection of the spherical image along $\ee{3}$, where $\ff{1} \p{t}$ is extended with a dashed line to the great circle it traces out.  Note that the great circle containing the poles and $\ee{3}$ is at an angle $\tau$ from the great circle containing $\ee{2}$, and $\ee{3}$, because $\ell$ is at a right angle to $\ff{1}$.}
   \label{fig:frame-path}
\end{figure}

\prf{%
  Following from \eqref{eq:ff1}, the point $\ff{1} \p{t}$ lies on the great circle of $\ee{1} \p{t}$ and $\ee{2} \p{t}$.  In particular, the arc length between $\ee{1} \p{t}$ and $\ff{1} \p{t}$ is the constant $\tau$ for all $t$.  And, from the Frenet-Serret equations, it follows that $\dot{\bm{e}}_{1} \p{t}$ is parallel to $\ee{2} \p{t}$, so $\ee{1} \p{t}$ must, for all $t$, be tangent to the great-circular arc between $\ee{1} \p{t}$ and $\ff{1} \p{t}$.  Putting those facts together, it follows that $\ee{1}$ is a spherical tractrix along $\ff{1}$; see \ref{fig:frame-path}, center.

  Similarly, $\dot{\bm{e}}_{3} \p{t}$ is parallel to $\ee{2} \p{t}$, making the great circle containing $\ee{2} \p{t}$, and $\ee{3} \p{t}$, tangent to $\ee{3} \p{t}$ for all $t$.  Also, note that $\ff{1} \p{t}$ traces out the great circle that is the equator determined by the axis $\ell$, because of its definition as the tangent to the parallels of $\Sigma$.  From those facts, along with looking along the direction of $\ee{3}$, it is quick to see that $\ee{3} \p{t}$ intersects the meridians (relative to $\ell$) at the constant angle $\tau$, making $\ee{3} \p{t}$ a loxodrome; see \ref{fig:frame-path}, right.
}

Now, as $\ee{1} \p{t}$ is a spherical tractrix, it has a cusp; \ref{fig:frame-path} shows what this looks like for the catenoid.  And, at that cusp, the Frenet-Serret frame is not defined, and it results in a discontinuity of $\ee{2} \p{t}$, and thusly, $\ee{3} \p{t}$, whereby they get reflected about the origin.  In order to fix these discontinuities: without loss of generality, assume the cusp occurs at $t = 0$, and then, let
\begin{align*}
  \eeh{2} \p{t} &= \begin{cases}
    - \ee{2} \p{t} &t < 0\\
    \displaystyle \lim_{t \to 0^{-}} - \ee{2} \p{t} = \lim_{t \to 0^{+}} \ee{2} \p{t} &t = 0\\
    \ee{2} \p{t} &t > 0
    \end{cases}\\
  \eeh{3} \p{t} &= \ee{1} \p{t} \times \eeh{2} \p{t} \mpn{,}
\end{align*}
wherein the limits exist because there is a removable discontinuity at $t = 0$.  As such, there is now a continuous framing $\cur{\ee{1}, \eeh{2}, \eeh{3}}$ of the curve $\alpha$, where $\eeh{3} \p{t}$ is still a loxodrome, as a reflection does not change that property.

Next, in order to find a differential equation for $\alpha$, the curves $\ee{1}$ and $\eeh{3}$ will first be found in the following lemma.

\lemma{%
  The curves $\ee{1}, \eeh{3} \colon \R \fdtr \Sph^{2}$ are given by
  \begin{align*}
    \ee{1} \p{t} &= \p{- \cos t \sin \tau \tanh \p{t \cot \tau} + \sin t \cos \tau,\right.\\
      &\hspace{15ex}\left. - \sin t \sin \tau \tanh \p{t \cot \tau} - \cos t \cos \tau ,\msp \sin \tau \sech \p{t \cot \tau}}\\
    \eeh{3} \p{t} &= \f{1}{\cosh \p{t \cot \tau}} \p{\cos t,\msp \sin t,\msp \sinh \p{t \cot \tau}} \mpn{.}
  \end{align*}
}

\prfeq{%
  First $\eeh{3}$ will be determined, from which $\ee{1}$ will then be determined.  For ease of computation, consider the stereographic projection $\eec{3} \p{t} \defeq \sigma \circ \eeh{3} \p{t}$ of $\eeh{3} \p{t}$, wherein
  \begin{align*}
    \sigma \colon \Sph^{2} &\fdtr \C\\
      \p{x, y, z} &\Mapsto \f{x}{1 - z} \;+\; \im \f{y}{1 - z} \mpn{.}
  \end{align*}
  Because $\eeh{3} \p{t}$ is a loxodrome, it follows from the properties of the stereographic projection $\sigma$, that $\eec{3} \p{t}$ meets the pencil of lines through the origin at the constant angle $\tau$: namely,
  \delimitershortfall=5pt
  \[\f{\dot{\check{\bm{e}}}_{3} \ol{\check{\bm{e}}}_{3}}{\abs{\dot{\check{\bm{e}}}_{3}} \abs{\vphantom{\bar{\check{\bm{e}}}_{3}}\smash{\ol{\check{\bm{e}}}_{3}}}} \overset{!}{=} \e^{\im \tau} \mpn{.}\]
  \delimitershortfall=0pt
  After letting $\eec{3} \p{t} = r \p{t} \e^{\im \theta \p{t}}$, that is
  \[\f{\dot{r} \e^{\im \theta} + \im r \e^{\im \theta} \dot{\theta}}{\sqrt{\dot{r}^{2} + r^{2} \dot{\theta}^{2}}} \; \f{r \e^{-\im \theta}}{r} \overset{!}{=} \e^{\im \tau} \mpn{,}\]
  which implies
  \[\f{\dot{r}}{\sqrt{\dot{r}^{2} + r^{2} \dot{\theta}^{2}}} \overset{!}{=} \cos \tau \quad\text{and}\quad \f{r \dot{\theta}}{\sqrt{\dot{r}^{2} + r^{2} \dot{\theta}^{2}}} \overset{!}{=} \sin \tau \mpn{,}\]
  that together make
  \[\f{\dot{r}}{r} \overset{!}{=} \dot{\theta} \cot \tau \mpn{.}\]
  Thus,
  \[r \p{t} = A \e^{\theta \p{t} \cot \tau} \quad\text{so}\quad \eec{3} \p{t} = A \e^{\p{\cot \tau + \im} \theta \p{t}} \mpn{,}\]
  for some constant $A$ and some parametrization $\theta \colon \R \fdtr \R$, which will be taken to be the identity for simplicity.  With the initial condition that $\alpha \p{0}$ is at the neck of $\Sigma$, it follows that $\eeh{3} \p{0}$ lies on the great circle perpendicular to $\ell$, and thus, that $\eec{3} \p{0}$ should lie on the unit circle of $\C$.  Therewith, it is sufficient to set $A = 1$.  Finally,
  \begin{align*}
    \eeh{3} \p{t} &= \sigma^{-1} \circ \eec{3} \p{t}\\
      &= \f{1}{\e^{2 t \cot \tau} + 1} \p{2 \e^{t \cot \tau} \cos t,\msp 2 \e^{t \cot \tau} \sin t,\msp \e^{2 t \cot \tau} - 1}\\
      &= \f{1}{\cosh \p{t \cot \tau}} \p{\cos t,\msp \sin t,\msp \sinh \p{t \cot \tau}} \mpn{.}
  \end{align*}

  Now, with that, $\ee{1}$ will determined.  To that end, first note that the Frenet-Serret equations show that $\dot{\bm{e}}_{3} \p{t} \| \ee{2} \p{t}$ for all $t$ outside of the cusp on $\ee{1}$, and thusly, that $\dot{\hat{\bm{e}}}_{3} \p{t} \| \eeh{2}$ for all $t$.  Then, as such, it follows that the span of $\ee{1} \p{t}$ lies in the envelope $\Omega$ of planes normal to $\eeh{3} \p{t}$.  The envelope $\Omega$ is obtained by intersecting the plane normal to $\eeh{3} \p{t}$,
  \[\Pi_{\perp} \p{t} = \cur{\bm{x} = \p{x^{1}, x^{2}, x^{3}} \in \R^{3} \colon\; \subnote{\cosh \p{t \cot \tau} \eeh{3} \p{t} \cdot \bm{x}}{50ex}{$x^{1} \cos t + x^{2} \sin t + x^{3} \sinh \p{t \cot \tau}$} = 0}\ubfixvsp \mpn{,}\]
  with its derivative,
  \[\dot{\Pi}_{\perp} \p{t} = \cur{\bm{x} \in \R^{3} \colon\; - x^{1} \sin t + x^{2} \cos t + x^{3} \cot \tau \cosh \p{t \cot \tau} = 0} \mpn{,}\]
  which will always be a line because $\dot{\hat{\bm{e}}}_{3} \| \eeh{2}$.

  To determine a parametrization of $\ee{1}$, it suffices to take the intersection of $\Omega$ with any plane, that is nowhere tangent to $\Omega$, and then, to project that intersection onto $\Sph^{2} \subset \R^{3}$.  So, taking the plane $\Pi_{1} \subset \R^{3}$ such that $x^{3} = 1$, suffices as it is tangent to $\Sph^{2}$, and thus, cannot be anywhere tangent to $\Omega$.  Doing so, the intersection $\ees{1} \p{t} \defeq \Pi_{\perp} \p{t} \cap \dot{\Pi}_{\perp} \p{t} \cap \Pi_{1}$ is: $\bm{x} \in \R^{3}$ such that
  \begin{align*}
    &\left\{\begin{aligned}
    x^{1} \cos t + x^{2} \sin t + \sinh \p{t \cot \tau} &= 0\\
    - x^{1} \sin t + x^{2} \cos t + \cot \tau \cosh \p{t \cot \tau} &= 0\\
    x^{3} &= 1
    \end{aligned}\right.\\
    \iff &\left\{\begin{aligned}
    x^{1} &= - \cos t \sinh \p{t \cot \tau} + \sin t \cot \tau \cosh \p{t \cot \tau}\\
    x^{2} &= - \sin t \sinh \p{t \cot \tau} - \cos t \cot \tau \cosh \p{t \cot \tau}\\
    x^{3} &= 1
    \end{aligned}\right. \mpn{.}
  \end{align*}
  Normalizing, obtains on $\Sph^{2}$
  \reqnos
  \begin{align*}
    \ee{1} \p{t} &= \f{\ees{1} \p{t}}{\abs{\ees{1} \p{t}}}\\
      &\hspace{-5ex}= \f{1}{\csc \tau \cosh \p{t \cot \tau}} \p{- \cos t \sinh \p{t \cot \tau} + \sin t \cot \tau \cosh \p{t \cot \tau} ,\right.\\
        &\hphantom{\f{1}{\csc \tau \cosh \p{t \cot \tau}}}\hspace{2ex}\left. - \sin t \sinh \p{t \cot \tau} - \cos t \cot \tau \cosh \p{t \cot \tau} , \msp 1}\\
      &\hspace{-5ex}= \p{- \cos t \sin \tau \tanh \p{t \cot \tau} + \sin t \cos \tau,\right.\\
        &\hspace{5ex}\left. - \sin t \sin \tau \tanh \p{t \cot \tau} - \cos t \cos \tau ,\msp \sin \tau \sech \p{t \cot \tau}} \mpn{.}\eqprf
  \end{align*}
  \leqnos
}

Now, using that lemma, the following proposition finds the asymptotic curve $\alpha$.

\prop{%
  Assume that the asymptotic curve $\alpha \colon \R \fdtr \R^{3}$ has the form
  \[\alpha \p{t} = \p{r \p{t} \cos t ,\msp r \p{t} \sin t ,\msp h \p{t}} \mpn{,}\]
  for some functions $r, h \colon \R \fdtr \R$, corresponding to the radius, and height, functions, respectively.  Then,
  \begin{align*}
    r \p{t} &= J \cosh^{\tan^{2} \tau} \p{t \cot \tau}\\
    h \p{t} &= - J \tan \tau \int \cosh^{- 1 + \tan^{2} \tau} \p{t \cot \tau} \;\di t \mpn{,}
  \end{align*}
  for some $J \in \R$, which is the same the result gotten in \eqref{eq:profile_m} but after a reparametrization.
}

\prf{%
  It is possible to find $\alpha \p{t}$ from $\ees{1} \p{t}$: so that $\alpha$ is an asymptotic curve, it suffices that its tangents lie along the line spanned by $\ee{1}$,
  \begin{align*}
    \p{\dot{r} \p{t} \cos t - r \p{t} \sin t ,\msp \dot{r} \p{t} \sin t + r \p{t} \cos t ,\msp \dot{h} \p{t}} &= \dot{\alpha} \p{t}\\
    &\overset{!}{=} \lambda \p{t} \ees{1} \p{t} \mpn{,}
  \end{align*}
  wherein $\lambda \colon \R \fdtr \R \setminus \cur{0}$.  This is the following system of equations
  \[\mathclap{\left\{\begin{aligned}
      \dot{r} \p{t} \cos t - r \p{t} \sin t &= \lambda \p{t} \p{- \cos t \sinh \p{t \cot \tau} + \sin t \cot \tau \cosh \p{t \cot \tau}}\\
      \dot{r} \p{t} \sin t + r \p{t} \cos t &= \lambda \p{t} \p{- \sin t \sinh \p{t \cot \tau} - \cos t \cot \tau \cosh \p{t \cot \tau} }\\
      \dot{h} \p{t} &= \lambda \p{t}
    \end{aligned}\right.}\]
  \begin{align*}
    \iff &\left\{\begin{aligned}
      r \p{t} &= - \lambda \p{t} \cot \tau \cosh \p{t \cot \tau}\\
      \dot{r} \p{t} &= - \lambda \p{t} \sinh \p{t \cot \tau}\\
      \dot{h} \p{t} &= \lambda \p{t}
    \end{aligned}\right.\\
    \iff &\left\{\begin{aligned}
      \f{\dot{r} \p{t}}{r \p{t}} &= \tan \tau \tanh \p{t \cot \tau}\\
      \dot{h} \p{t} &= - \f{r \p{t}}{\cot \tau \cosh \p{t \cot \tau}}
    \end{aligned}\right.\\
    \iff &\left\{\begin{aligned}
      \log r \p{t} &= \tan^{2} \tau \log \p{\cosh \p{t \cot \tau}} + \log J &&\hspace{1.3ex}\mnote{20ex}{for some constant $J > 0$}\\
      \dot{h} \p{t} &= - \f{r \p{t}}{\cot \tau \cosh \p{t \cot \tau}}
    \end{aligned}\right.\\
    \iff &\left\{\begin{aligned}
      r \p{t} &= J \cosh^{\tan^{2} \tau} \p{t \cot \tau} &&\mnote{20ex}{for some constant $J > 0$}\\
      h \p{t} &= - J \tan \tau \int \cosh^{- 1 + \tan^{2} \tau} \p{t \cot \tau} \;\di t
    \end{aligned}\right. \hspace{-3ex}\!\phantom{,}\mpn{,}
  \end{align*}
  where the constant $J$ is simply the scaling of $\alpha$.  Lastly, reparametrizing with
  \[t = \sgn \theta \ \tan \tau \ \arccosh \sec \theta \mpn{,}\]
  and recalling that $m = \tan^{2} \tau$, obtains \eqref{eq:profile_m} with $c = 0 = K'$.
}

And therewith, it is possible to obtain the asymptotic parametrization of $\Sigma$.

\prop{%
  An asymptotic parametrization of $\p{m , 0}$-type is given by
  \begin{align*}
    \R \times \R &\fdtr \R^{3}\\
    \p{t , s} &\Mapsto \p{r \p{t - s} \cos \p{t + s} ,\msp r \p{t - s} \sin \p{t + s} ,\msp h \p{t -s}} \mpn{,}
  \end{align*}
  wherein $r , h : \R \fdtr \R$ are given in \ref{5.3}.
}

\begin{figure}[h!]
  \centering{%
    \begin{minipage}[m]{0.45\textwidth}
      \includegraphics[width=\textwidth]{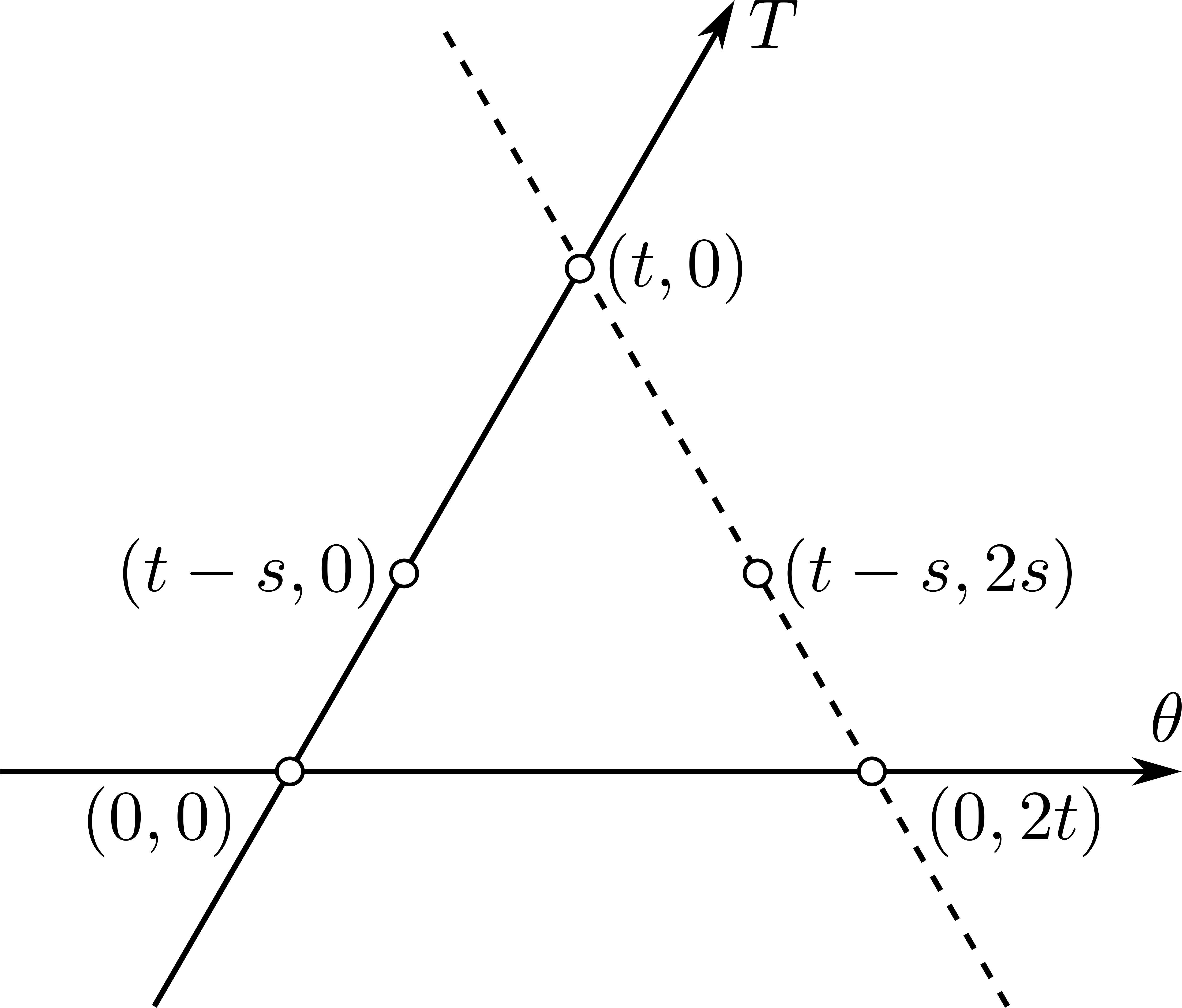}
    \end{minipage}
  }
  \caption{%
    In the $\p{t , \theta}$-coordinate space, the dotted line corresponds to an asymptotic curve in the other family.  To trace along that dotted line depending on $s$, the coordinates would be $\p{t - s , 2 s}$ for a given $t$.
  }
  \label{fig:asymptotic-param}
\end{figure}

\prfeq{%
  From the symmetry of $\p{m , 0}$-type being a surface of revolution, it quickly follows that the other family of asymptotic curves are gotten by a reflection in $\R^{3}$.  In particular, this means that they can be simply gotten by a reparametrization of the family of asymptotic curves from \ref{5.3}: let the parameter of the family be $\theta \in \R$; then, the curves in the family are
  \[\p{r \p{T} \cos \p{T + \theta} ,\msp r \p{T} \sin \p{T + \theta} ,\msp h \p{T}} \mpn{,}\]
  where $T \in \R$ is the parameter along the curves.  Using \ref{fig:asymptotic-param}, an appropriate reparametrization would be
  \begin{align*}
    T \p{t , s} &= t - s\\
    \theta \p{t , s} &= 2 s \mpn{,}
  \end{align*}
  which then makes the asymptotic parametrization
  \reqnos
  \begin{equation*}
    \p{r \p{t - s} \cos \p{t + s} ,\msp r \p{t - s} \sin \p{t + s} ,\msp h \p{t - s}} \mpn{.}\eqprf
  \end{equation*}
  \leqnos
}

Finally, the following remark closes the section, discussing how a constant-angle parametrization for surfaces with $m < 0$ was simultaneously obtained.

\remark{%
As mentioned in \ref{5.3}, the results obtained there (for $m > 0$) are, after reparametrization, the same as the results in \eqref{eq:profile_m}, which hold for any sign of $m$.  It is therewith suggested that the assumption that the surface is negatively-curved ($m > 0$) was somehow immaterial.  In the case that $m < 0$, the objective is to find a curve $\alpha$ which meets the parallel circles of $\Sigma$ at a constant angle $\tau$.  But, as $\Sigma$ is positively-curved ($m < 0$), this curve cannot be an asymptotic curve, so its osculating planes do not agree with the tangent planes of $\Sigma$ and the unit tangent does not trace out a spherical tractrix.  However, it is still the case that the unit normal along $\alpha$ traces out a loxodrome, as in \ref{5.2}.  That loxodrome can be used just as in \ref{5.3}, and then \ref{5.4}, to obtain a constant-angle parametrization for $\Sigma$ in this case.  Correspondingly, there is an imaginary $\tau$ so that $m = \tan^{2} \tau$: this can be justified using Laguerre's formula (see \cite[\S18.7]{RichterGebert}) on intersection of the projectivized Dupin indicatrix conic with the line at infinity.}

\section{Explicit examples for $m \in \Z$}
\label{sec:examples}

In this section, there will be examples of profile curves from previous section's four cases (even/odd negative/positive values of $m$).

\Exx{%
  From \ref{4.1}: for $m = 2$,
  \begin{align*}
    r \p{\theta} &= J \sec^{2} \theta + \f{c}{3} \cos \theta\\
    h \p{\theta} &= - 2 J \tan \theta + \f{c}{3} \sin \theta + K \mpn{.}
  \end{align*}
  The satisfies the following algebraic relationship, when $K = 0$,
  \begin{align*}
    &c^{6} - 27 c^{4} h^{2} + 243 c^{2} h^{4} - 729 h^{6} + 135 c^{4} J^{2} + 1944 c^{2} h^{2} J^{2} - 8748 h^{4} J^{2} + 3888 c^{2} J^{4}\\
    &\hspace{5ex}- 34992 h^{2} J^{4} - 46656 J^{6} - 54 c^{4} J r - 486 c^{2} h^{2} J r + 8748 h^{4} J r - 1944 c^{2} J^{3} r\\
    &\hspace{5ex}+ 69984 h^{2} J^{3} r + 139968 J^{5} r - 9 c^{4} r^{2} + 162 c^{2} h^{2} r^{2} - 729 h^{4} r^{2} - 1296 c^{2} J^{2} r^{2}\\
    &\hspace{5ex}- 40824 h^{2} J^{2} r^{2} - 151632 J^{4} r^{2} + 648 c^{2} J r^{3} + 5832 h^{2} J r^{3} + 69984 J^{3} r^{3}\\
    &\hspace{5ex}- 11664 J^{2} r^{4} = 0 \mpn{;}
  \end{align*}
  when $c = 0 = K$, this can be reduced to
  \[4 J r - h^{2} - 4 J^{2} = 0 \mpn{,}\]
  showing that profile curve is a parabola.  Its evolute is
  \[\epsilon \p{\theta} = \p{3 J \sec^{2} \theta ,\msp 2 J \tan^{3} \theta} \mpn{,}\]
  which is also algebraic, satisfying
  \delimitershortfall=5pt
  \[27 \p{\epsilon^{2}}^{2} J + 108 J^{3} - 108 J^{2} \p{\epsilon^{1}} + 36 J \p{\epsilon^{1}}^{2} - 4 \p{\epsilon^{1}}^{3} = 0 \mpn{.}\]
  \delimitershortfall=0pt

  As is the case for the other such $m$, as $\theta \to \pm \f{\pi}{2}$, the coordinate functions $r, h \to \pm \infty$; furthermore, translating by $\pi$ in the domain, $r \p{\theta + \pi}$ and $h \p{\theta + \pi}$, is the same as flipping the sign of $c$.  Thus, the image of these profile curves on $\R$ is invariant under flipping the sign of $c$.\\[1ex]
  \begin{minipage}[]{0.5\textwidth}
    \includegraphics[width=\textwidth]{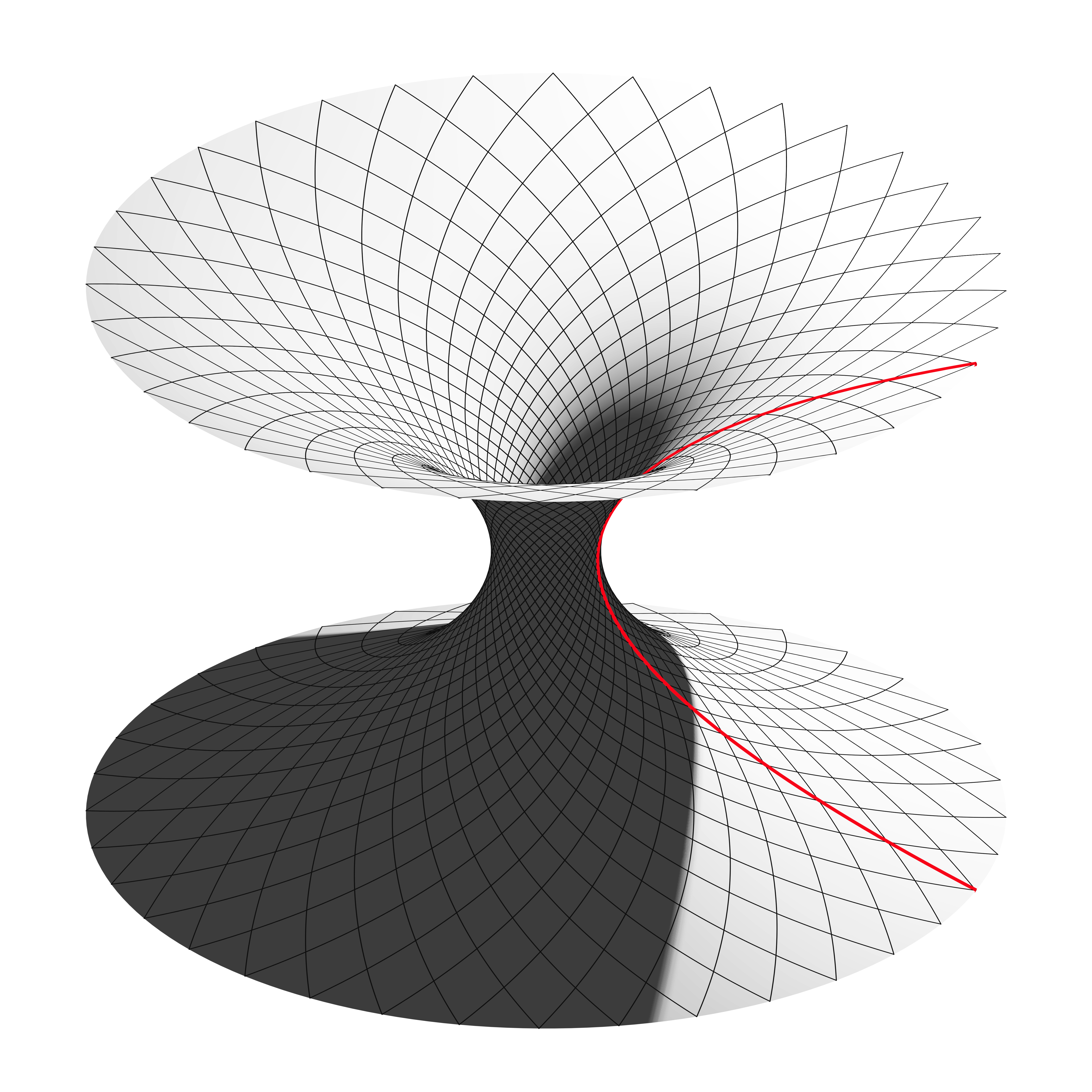}
  \end{minipage}
  \hfill
  \begin{minipage}[]{0.5\textwidth}
    \includegraphics[width=\textwidth]{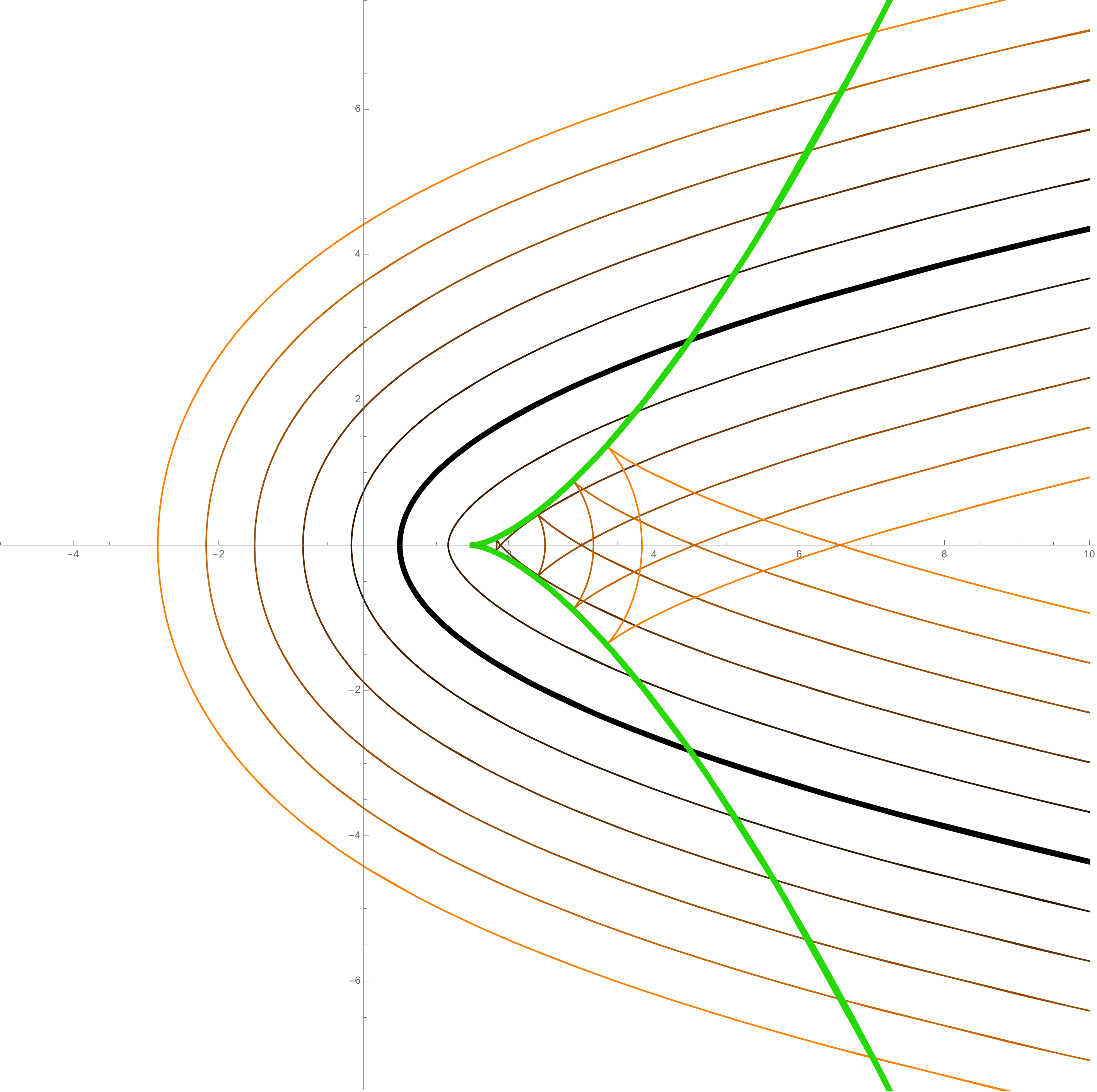}
  \end{minipage}\\[1ex]
  On the left, there is a rendering of a surface of $\p{2 , 0}$-type in its asymptotic parametrization from \Aref{sec:param}, with its parabolic profile curve shown in red.  On the right, the profile curves when $J = \f{1}{2}$ the profile curves for various $c$ are shown, with the curve with $c = 0$ in bold black, and their common evolute in green.  Note that, for a given $c \ne 0$, the profile curve has two components.
}

\Exx{%
  From \ref{4.2}: for $m = -3$,
  \begin{align*}
    r \p{\theta} &= J \cos^{3} \theta - \f{c}{2} \cos \theta\\
    h \p{\theta} &= 3 J \p{\sin \theta - \f{\sin^{3} \theta}{3}} - \f{c}{2} \sin \theta + K \mpn{.}
  \end{align*}
  The satisfies the following algebraic relationship, when $K = 0$,
  \begin{align*}
    &c^{2} h^{4} - 4 h^{6} + 4 c^{3} h^{2} J - 24 c h^{4} J + 4 c^{4} J^{2} - 56 c^{2} h^{2} J^{2} + 48 h^{4} J^{2} - 48 c^{3} J^{3}\\
    &\hspace{1cm}+ 192 c h^{2} J^{3} + 208 c^{2} J^{4} - 192 h^{2} J^{4} - 384 c J^{5} + 256 J^{6} + 2 c^{2} h^{2} r^{2}\\
    &\hspace{1cm}- 12 h^{4} r^{2} - 4 c^{3} J r^{2} - 12 c h^{2} J r^{2} + 16 c^{2} J^{2} r^{2} + 96 h^{2} J^{2} r^{2} + 48 c J^{3} r^{2}\\
    &\hspace{1cm}- 192 J^{4} r^{2} + c^{2} r^{4} - 12 h^{2} r^{4} + 12 c J r^{4} - 60 J^{2} r^{4} - 4 r^{6} = 0 \mpn{;}
  \end{align*}
  when $c = 0 = K$, this reduces to
  \begin{multline*}
    h^{6} - 12 h^{4} J^{2} + 48 h^{2} J^{4} - 64 J^{6} + 3 h^{4} r^{2}\\
      - 24 h^{2} J^{2} r^{2} + 48 J^{4} r^{2} + 3 h^{2} r^{4} + 15 J^{2} r^{4} + r^{6} = 0 \mpn{.}\!\phantom{.}
  \end{multline*}
  Its evolute is
  \[\epsilon \p{\theta} = \p{- 2 J \cos^{3} \theta ,\msp 2 J \sin^{3} \theta} \mpn{,}\]
  which is also algebraic, satisfying
  \delimitershortfall=5pt
  \begin{align*}
    &\p{\epsilon^{2}}^{6} - 12 \p{\epsilon^{2}}^{4} J^{2} + 48 \p{\epsilon^{2}}^{2} J^{4} - 64 J^{6} + 3 \p{\epsilon^{2}}^{4} \p{\epsilon^{1}}^{2} + 84 \p{\epsilon^{2}}^{2} J^{2} \p{\epsilon^{1}}^{2}\\
    &\hspace{1cm}+ 48 J^{4} \p{\epsilon^{1}}^{2} + 3 \p{\epsilon^{2}}^{2} \p{\epsilon^{1}}^{4} - 12 J^{2} \p{\epsilon^{1}}^{4} + \p{\epsilon^{1}}^{6} = 0 \mpn{.}
  \end{align*}
  \delimitershortfall=0pt

  As is the case for the other such $m$, as $\theta \to \pm \f{\pi}{2}$, the coordinate functions $r, h \to \pm \infty$; furthermore, translating by $\pi$ in the domain, $r \p{\theta + \pi}$ and $h \p{\theta + \pi}$, is the same as flipping the sign of $c$.  Thus, the image of these profile curves on $\R$ is invariant under flipping the sign of $c$.\\[1ex]
  \begin{minipage}[]{0.5\textwidth}
    \includegraphics[width=\textwidth]{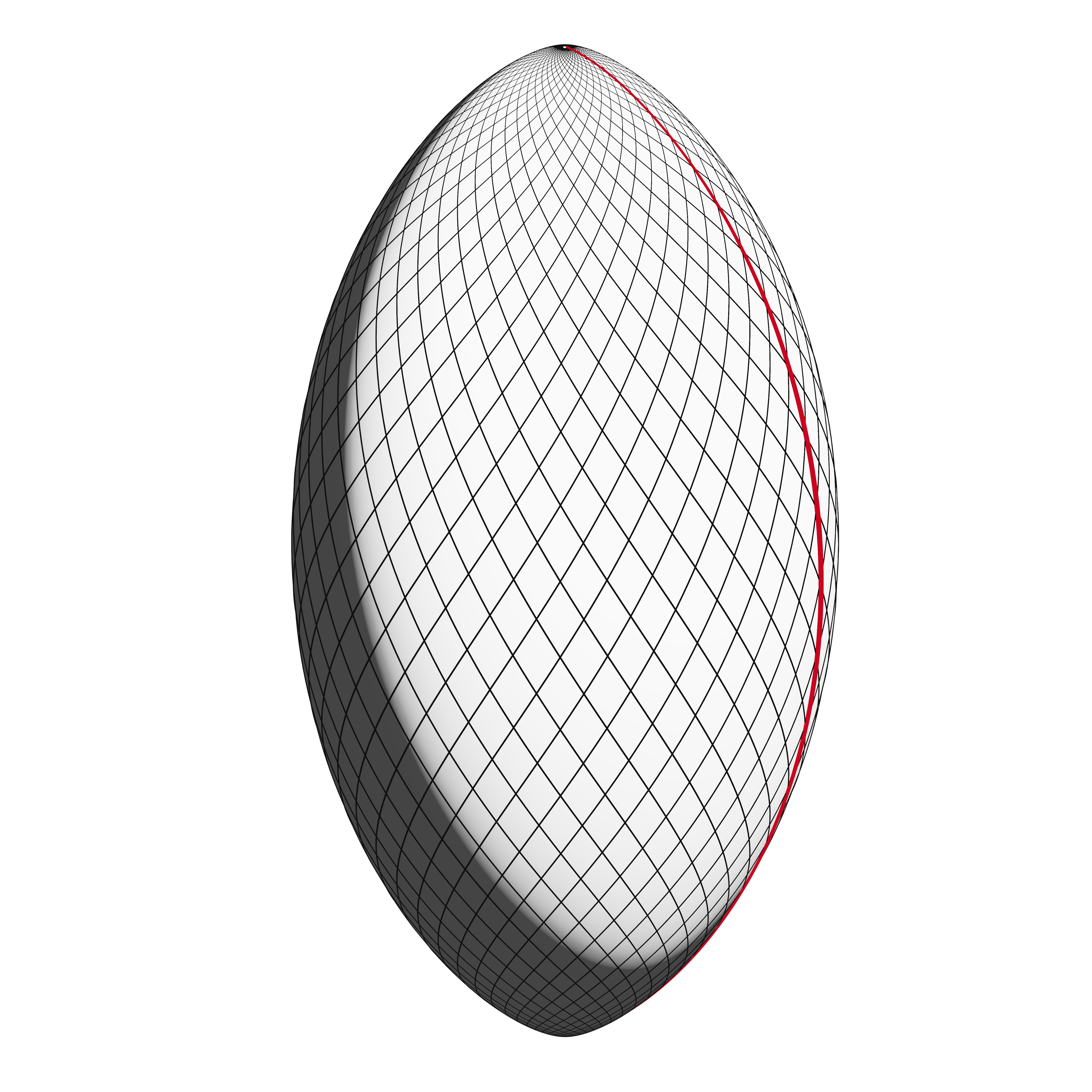}
  \end{minipage}
  \begin{minipage}[]{0.5\textwidth}
    \includegraphics[width=\textwidth]{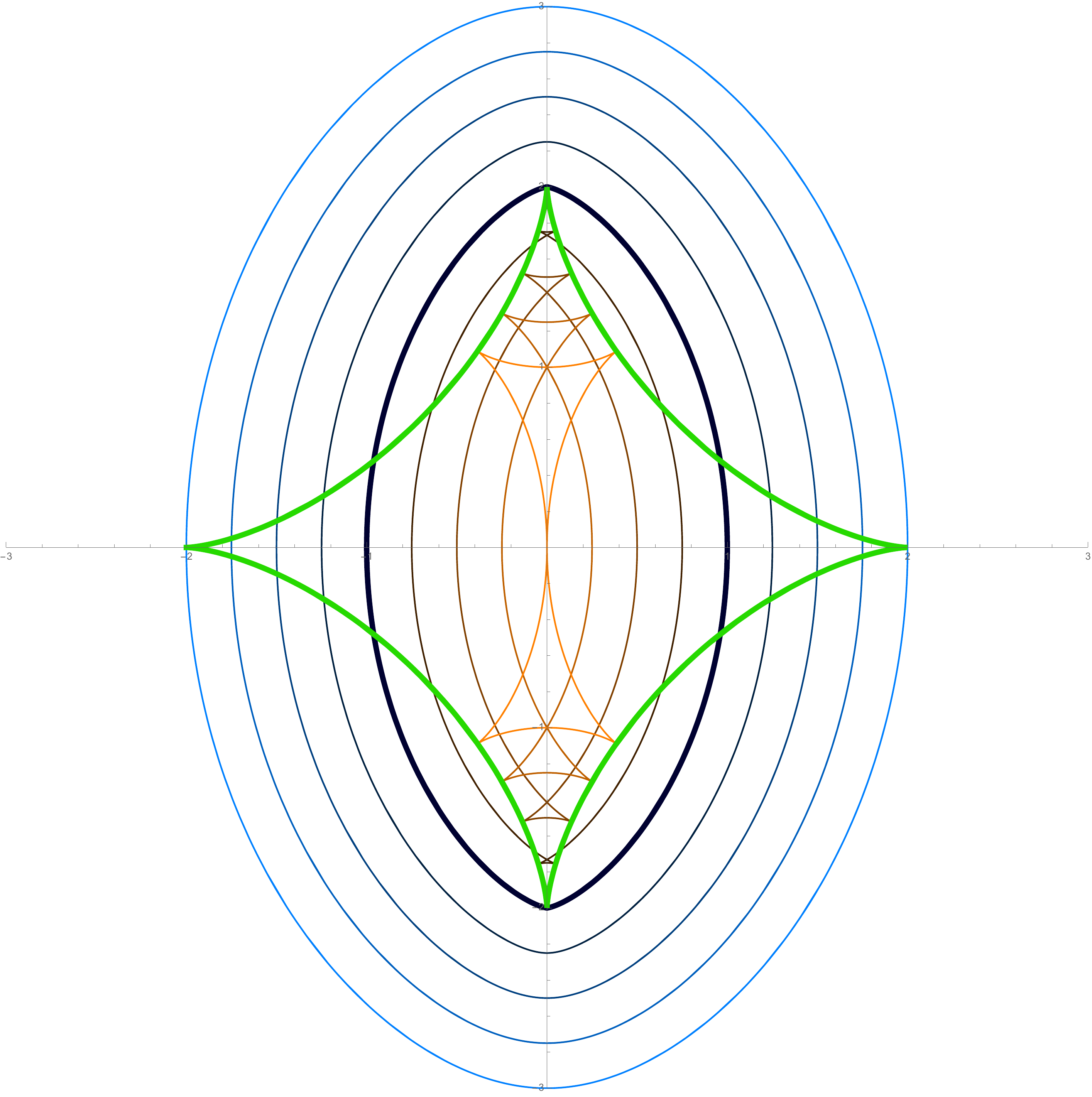}
  \end{minipage}\\[1ex]
  On the left, there is a rendering of a surface of $\p{-3 , 0}$-type with the constant-angle parametrization as discussed in \ref{5.5}, with its profile curve shown in red.  On the right, the profile curves when $J = 1$ for various $c$ are shown, with the curve with $c = 0$ in bold black, and their common evolute in green.  As $c$ decreases, the profile curves comes more blue; while as $c$ increase, they become more orange.
}

\Exx{%
  From \ref{4.3}: for $m = 3$,
  \begin{align*}
    r \p{\theta} &= J \sec^{3} \theta + \f{c}{4} \cos \theta\\
    h \p{\theta} &= - \f{3 J}{2} \p{\sec \theta \tan \theta + \log \abs{\sec \theta + \tan \theta}} + \f{c}{4} \sin \theta + K \mpn{.}
  \end{align*}

  As is the case for the other such $m$, as $\theta \to \pm \f{\pi}{2}$, the coordinate functions $r, h \to \pm \infty$; furthermore, translating by $\pi$ in the domain, $r \p{\theta + \pi}$ and $h \p{\theta + \pi}$, is the same as flipping the sign of $r$ and $h$.

  \begin{minipage}[]{0.475\textwidth}
    Here, for $J = \f{1}{2}$, the profile curves for various $c$ are shown, with the curve with $c = 0$ in bold black, and their common evolute in green.  As~$c$~\mbox{decreases}, the profile curves comes more blue; while as $c$ \mbox{increase}, they become more orange.  Note that the profile curves have two \mbox{components}.
  \end{minipage}
  \hfill
  \begin{minipage}[]{0.45\textwidth}
    \includegraphics[width=\textwidth]{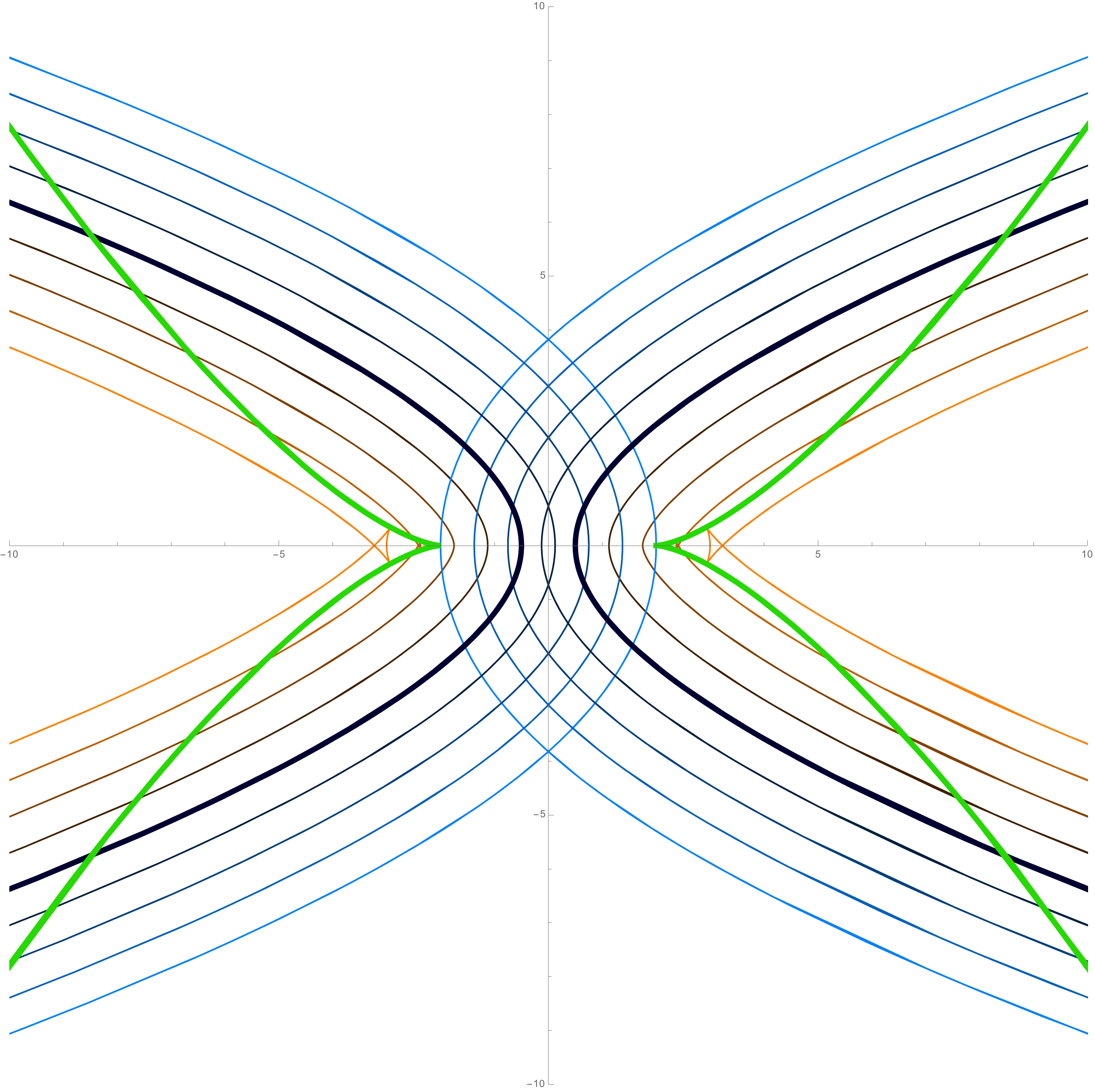}
  \end{minipage}\\
  \vspace{-3.5ex}
}

\Exx{%
  From \ref{4.4}: for $m = -2$,
  \begin{align*}
    r \p{\theta} &= J \cos^{2} \theta - c \cos \theta\\
    h \p{\theta} &= J \p{\theta + \cos \theta \sin \theta} - c \sin \theta + K \mpn{.}
  \end{align*}

  As is the case for the other such $m$, translating by $\pi$ in the domain, $r \p{\theta + \pi}$ and $h \p{\theta + \pi}$, is the same as flipping the sign of $c$ and translating $h$ by $J \pi$.

  It can be shown that the profile curve when $c = 0$ is a cycloid, as a roulette of a circle of radius $J$ along the axis of revolution, and that its evolute is similarly so but along the lines $r = - J$.

  Here, for $J = 1$, the profile curves for various $c$ are shown, with the curve with $c = 0$ in bold black, and their common evolute in green.  As $c$ decreases, the profile curves comes more blue; while as $c$ increase, they become more orange.\\[1ex]
  \begin{center}
    \includegraphics[width=0.45\textwidth]{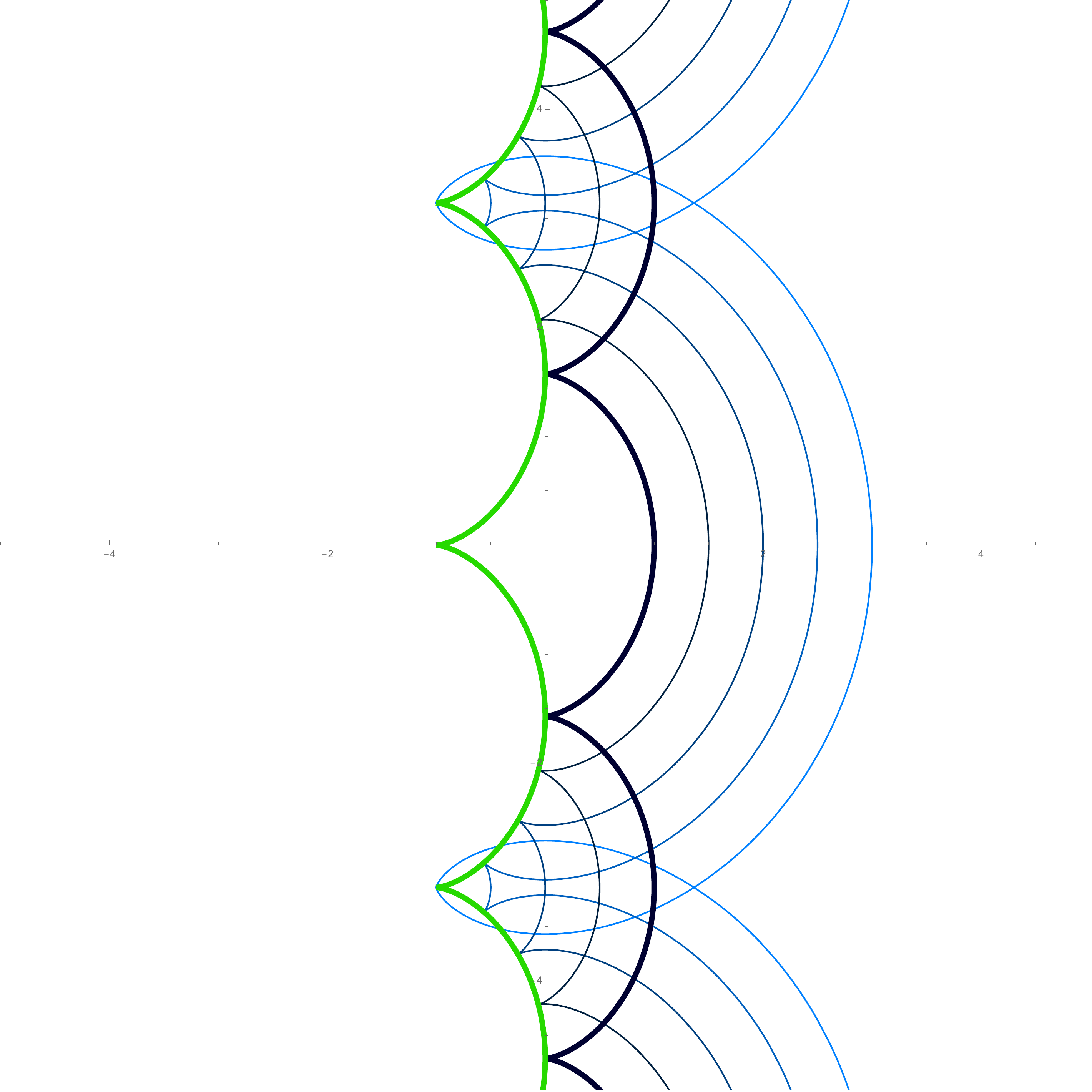}
    \hfill
    \includegraphics[width=0.45\textwidth]{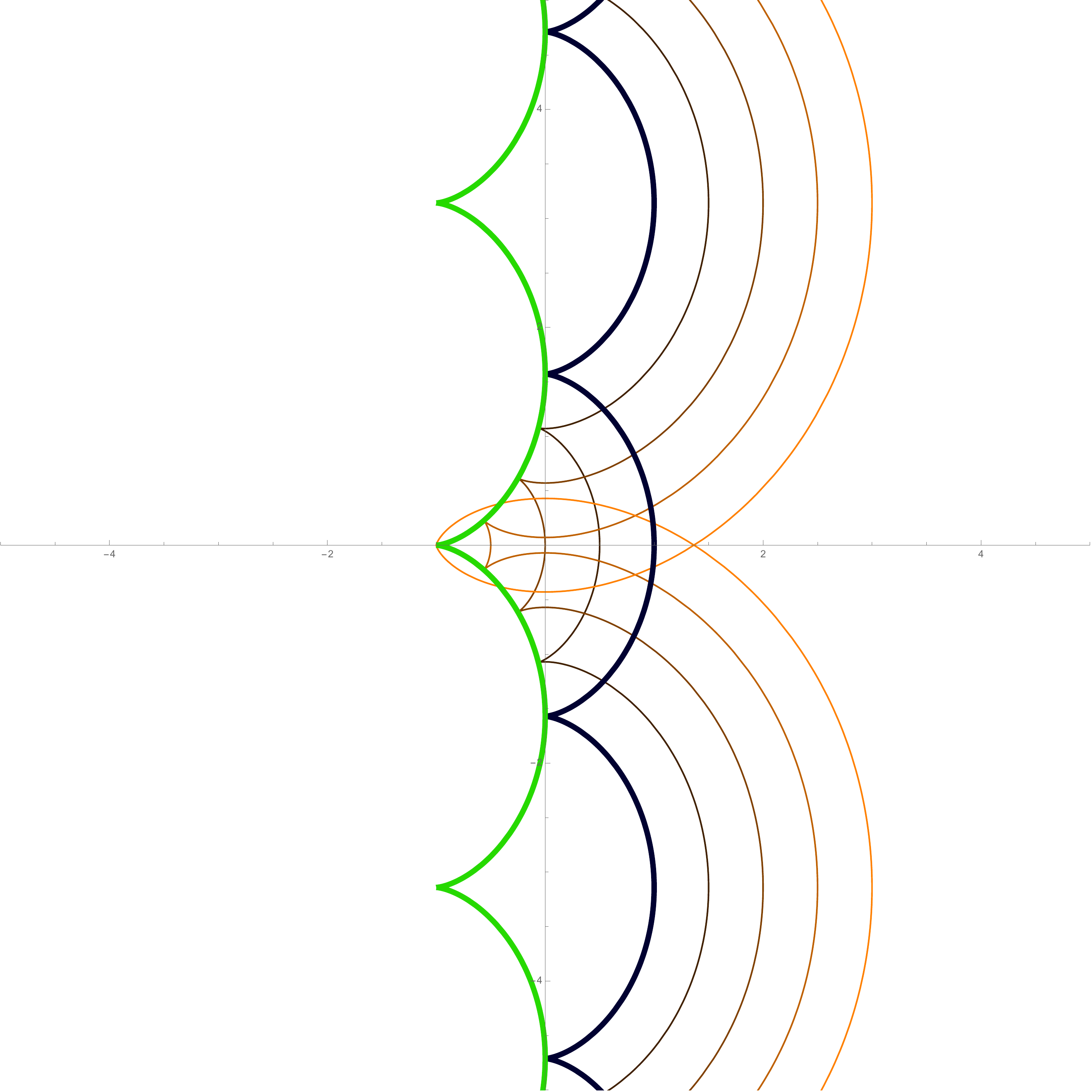}
  \end{center}
  \vspace{-3.5ex}
}

\section*{Acknowledgements}

The author is grateful to Helmut Pottmann for his notes that served as a foundation for this work.  To Christian M\"{u}ller, Arvin Rasoulzadeh, and Friedrich Manhart, the author is appreciative for the many insight conversations that were helpful in the development of this paper.  And, the author is thankful to Udo Hertrich-Jeromin for his feedback.

\section*{Funding}

This project partially received funding from the European Union's Horizon 2020 program under the Marie Sk\l{}odowska-Curie grant agreement No.~675789, and from the Austrian Science Fund (FWF) under the grant P~29981.

\biblio


\begin{bibdiv}
\begin{biblist}

\bib{Baikoussis-Koufogiorgos}{article}{
      author={Baikoussis, C.},
      author={Koufogiorgos, T.},
       title={On the inner curvature of the second fundamental form of
  helicoidal surfaces},
        date={1997},
        ISSN={0003-889X},
     journal={Arch. Math. (Basel)},
      volume={68},
      number={2},
       pages={169\ndash 176},
         url={https://doi.org/10.1007/s000130050046},
      review={\MR{1425508}},
}

\bib{Boyacioglu-Lopez_spacelike}{article}{
      author={Boyac\i~o\u{g}lu Kalkan, {\"{O}}zg{\"{u}}r},
      author={L{\'{o}}pez, Rafael},
       title={Spacelike surfaces in {M}inkowski space satisfying a linear
  relation between their principal curvatures},
        date={2011},
     journal={Differ. Geom. Dyn. Syst.},
      volume={13},
       pages={107\ndash 116},
      review={\MR{2812439}},
}

\bib{Boyacioglu-Lopez}{article}{
      author={Boyacioglu~Kalkan, {\"{O}}zg{\"{u}}r},
      author={L{\'{o}}pez, Rafael},
      author={Saglam, Derya},
       title={Linear {W}eingarten surfaces foliated by circles in {M}inkowski
  space},
        date={2011},
        ISSN={1027-5487},
     journal={Taiwanese J. Math.},
      volume={15},
      number={5},
       pages={1897\ndash 1917},
         url={https://doi.org/10.11650/twjm/1500406413},
      review={\MR{2880383}},
}

\bib{Ceyhan-Fokas-Guerses}{article}{
      author={Ceyhan, \"{O}.},
      author={Fokas, A.~S.},
      author={G\"{u}rses, M.},
       title={Deformations of surfaces associated with integrable
  {G}auss-{M}ainardi-{C}odazzi equations},
        date={2000},
        ISSN={0022-2488},
     journal={J. Math. Phys.},
      volume={41},
      number={4},
       pages={2251\ndash 2270},
         url={https://doi.org/10.1063/1.533237},
      review={\MR{1751887}},
}

\bib{Gray}{book}{
      author={Gray, Alfred},
      author={Abbena, Elsa},
      author={Salamon, Simon},
       title={Modern differential geometry of curves and surfaces with
  {M}athematica{$^\circledR$}},
     edition={Third},
      series={Studies in Advanced Mathematics},
   publisher={Chapman \& Hall/CRC, Boca Raton, FL},
        date={2006},
        ISBN={978-1-58488-448-4; 1-58488-448-7},
      review={\MR{2253203}},
}

\bib{Hadzhilazova-Mladenov}{article}{
      author={Hadzhilazova, Mariana},
      author={Mladenov, Iva{\"{i}}lo~M.},
       title={Once more the {M}ylar balloon},
        date={2008},
        ISSN={1310-1331},
     journal={C. R. Acad. Bulgare Sci.},
      volume={61},
      number={7},
       pages={847\ndash 856},
      review={\MR{2450943}},
}

\bib{Hopf}{article}{
      author={Hopf, Heinz},
       title={{\"{U}}ber {F}l{\"{a}}chen mit einer {R}elation zwischen den
  {H}auptkr{\"{u}}mmungen},
        date={1951},
        ISSN={0025-584X},
     journal={Math. Nachr.},
      volume={4},
       pages={232\ndash 249},
         url={https://doi.org/10.1002/mana.3210040122},
      review={\MR{40042}},
}

\bib{Jimenez-Mueller-Pottmann}{article}{
      author={Jimenez, Michael~R.},
      author={M\"{u}ller, Christian},
      author={Pottmann, Helmut},
       title={Discretizations of surfaces with constant ratio of principal
  curvatures},
        date={2020},
        ISSN={0179-5376},
     journal={Discrete Comput. Geom.},
      volume={63},
      number={3},
       pages={670\ndash 704},
         url={https://doi.org/10.1007/s00454-019-00098-7},
      review={\MR{4074339}},
}

\bib{Kuehnel}{article}{
      author={K{\"{u}}hnel, Wolfgang},
       title={Zur inneren {K}r{\"{u}}mmung der zweiten {G}rundform},
        date={1981},
        ISSN={0026-9255},
     journal={Monatsh. Math.},
      volume={91},
      number={3},
       pages={241\ndash 251},
         url={https://doi.org/10.1007/BF01301791},
      review={\MR{619967}},
}

\bib{KuehnelDiffGeom}{book}{
      author={K{\"{u}}hnel, Wolfgang},
       title={Differential geometry},
      series={Student Mathematical Library},
   publisher={American Mathematical Society, Providence, RI},
        date={2015},
      volume={77},
        ISBN={978-1-4704-2320-9},
         url={https://doi.org/10.1090/stml/077},
        note={Curves---surfaces---manifolds, Third edition [of MR1882174],
  Translated from the 2013 German edition by Bruce Hunt, with corrections and
  additions by the author},
      review={\MR{3443721}},
}

\bib{Kuehnel-Steller}{article}{
      author={K{\"{u}}hnel, Wolfgang},
      author={Steller, Michael},
       title={On closed {W}eingarten surfaces},
        date={2005},
        ISSN={0026-9255},
     journal={Monatsh. Math.},
      volume={146},
      number={2},
       pages={113\ndash 126},
         url={https://doi.org/10.1007/s00605-005-0313-4},
      review={\MR{2176338}},
}

\bib{Lopez_EandH}{article}{
      author={L{\'{o}}pez, Rafael},
       title={Linear {W}eingarten surfaces in {E}uclidean and hyperbolic
  space},
        date={2008},
        ISSN={0103-9059},
     journal={Mat. Contemp.},
      volume={35},
       pages={95\ndash 113},
      review={\MR{2584178}},
}

\bib{Lopez}{article}{
      author={L{\'{o}}pez, Rafael},
       title={On linear {W}eingarten surfaces},
        date={2008},
        ISSN={0129-167X},
     journal={Internat. J. Math.},
      volume={19},
      number={4},
       pages={439\ndash 448},
         url={https://doi.org/10.1142/S0129167X08004728},
      review={\MR{2416724}},
}

\bib{Lopez_Htype}{article}{
      author={L{\'{o}}pez, Rafael},
       title={Rotational linear {W}eingarten surfaces of hyperbolic type},
        date={2008},
        ISSN={0021-2172},
     journal={Israel J. Math.},
      volume={167},
       pages={283\ndash 301},
         url={https://doi.org/10.1007/s11856-008-1049-3},
      review={\MR{2448026}},
}

\bib{Lopez-Pampano}{article}{
      author={L{\'{o}}pez, Rafael},
      author={P{\'{a}}mpano, {\'{A}}lvaro},
       title={Classification of rotational surfaces in {E}uclidean space
  satisfying a linear relation between their principal curvatures},
        date={2020},
        ISSN={0025-584X},
     journal={Math. Nachr.},
      volume={293},
      number={4},
       pages={735\ndash 753},
         url={https://doi.org/10.1002/mana.201800235},
      review={\MR{4089078}},
}

\bib{Mladenov-Opera_balloon}{incollection}{
      author={Mladenov, Iva\"{\i}lo~M.},
      author={Oprea, John},
       title={The {M}ylar balloon: new viewpoints and generalizations},
        date={2007},
   booktitle={Geometry, integrability and quantization},
   publisher={Softex, Sofia},
       pages={246\ndash 263},
      review={\MR{2341209}},
}

\bib{Mladenov-Opera_deform}{incollection}{
      author={Mladenov, Iva\"{\i}lo~M.},
      author={Oprea, John},
       title={On some deformations of the mylar balloon},
        date={2007},
   booktitle={X{V} {I}nternational {W}orkshop on {G}eometry and {P}hysics},
      series={Publ. R. Soc. Mat. Esp.},
      volume={11},
   publisher={R. Soc. Mat. Esp., Madrid},
       pages={310\ndash 315},
      review={\MR{2504244}},
}

\bib{Pottmann}{unpublished}{
      author={Pottmann, Helmut},
       title={Rotational surfaces with a constant negative ratio of principal
  curvatures},
        date={May 2019},
        note={Unpublished manuscript, 4 pages},
}

\bib{Pulov-Hadzhilazova-Mladenov}{incollection}{
      author={Pulov, Vladimir~I.},
      author={Hadzhilazova, Mariana~Ts.},
      author={Mladenov, Iva\"{\i}lo~M.},
       title={The {M}ylar balloon: an alternative description},
        date={2015},
   booktitle={Geometry, integrability and quantization {XVI}},
   publisher={Avangard Prima, Sofia},
       pages={256\ndash 269},
         url={https://doi.org/10.7546/giq-16-2015-256-269},
      review={\MR{3363850}},
}

\bib{RichterGebert}{book}{
      author={Richter-Gebert, J{\"{u}}rgen},
       title={Perspectives on projective geometry},
   publisher={Springer, Heidelberg},
        date={2011},
        ISBN={978-3-642-17285-4},
         url={https://doi.org/10.1007/978-3-642-17286-1},
        note={A guided tour through real and complex geometry},
      review={\MR{2791970}},
}

\bib{Riveros-CorroI}{article}{
      author={Riveros, Carlos M.~C.},
      author={Corro, Armando M.~V.},
       title={Surfaces with constant {C}hebyshev angle},
        date={2012},
        ISSN={0387-3870},
     journal={Tokyo J. Math.},
      volume={35},
      number={2},
       pages={359\ndash 366},
         url={https://doi.org/10.3836/tjm/1358951324},
      review={\MR{3058712}},
}

\bib{Riveros-CorroII}{article}{
      author={Riveros, Carlos M.~C.},
      author={Corro, Armando M.~V.},
       title={Surfaces with constant {C}hebyshev angle {II}},
        date={2013},
        ISSN={0387-3870},
     journal={Tokyo J. Math.},
      volume={36},
      number={2},
       pages={379\ndash 386},
         url={https://doi.org/10.3836/tjm/1391177977},
      review={\MR{3161564}},
}

\bib{Staeckel}{misc}{
      author={St{\"{a}}ckel, Paul},
       title={{Beitr{\"{a}}ge zur Fl{\"{a}}chentheorie.}},
    language={German},
         how={{Leipz. Ber. 48, 478-504 (1896).}},
        date={1896},
}

\bib{Tang-Killian-Bo-Wallner-Pottmann}{incollection}{
      author={Tang, Chengcheng},
      author={Kilian, Martin},
      author={Bo, Pengbo},
      author={Wallner, Johannes},
      author={Pottmann, Helmut},
       title={Analysis and design of curved support structures},
        date={2016},
   booktitle={Advances in architectural geometry 2016},
      editor={Adriaenssens, Sigrid},
      editor={Gramazio, Fabio},
      editor={Kohler, Matthias},
      editor={Menges, Achim},
      editor={Pauly, Mark},
   publisher={VDF Hochschulverlag, ETH Z{\"u}rich},
       pages={8\ndash 23},
}

\end{biblist}
\end{bibdiv}
\end{document}